\documentclass[preprint,10pt]{elsarticle}

\usepackage{mathrsfs,amsmath,amssymb}
\usepackage{mathrsfs,amsmath,amssymb}
\usepackage{amscd} %可換図式の為に
\usepackage{relsize,amsmath} %文字サイズ変更の為に
\newcommand{\lsum}{\mathlarger{\sum}}

\makeatletter%\Bigg よりさらに大きいもの
\newcommand{\biggg}[1]{{\hbox{$\left#1\vbox to 20.5pt{}\right.\n@space$}}}
\newcommand{\Biggg}[1]{{\hbox{$\left#1\vbox to 23.5pt{}\right.\n@space$}}}
\newcommand{\bigggg}[1]{{\hbox{$\left#1\vbox to 26.5pt{}\right.\n@space$}}}
\newcommand{\Bigggg}[1]{{\hbox{$\left#1\vbox to 29.5pt{}\right.\n@space$}}}
\newcommand{\biggggg}[1]{{\hbox{$\left#1\vbox to 32.5pt{}\right.\n@space$}}}
\newcommand{\Biggggg}[1]{{\hbox{$\left#1\vbox to 35.5pt{}\right.\n@space$}}}
\newcommand{\bigggggg}[1]{{\hbox{$\left#1\vbox to 38.5pt{}\right.\n@space$}}}
\newcommand{\Bigggggg}[1]{{\hbox{$\left#1\vbox to 41.5pt{}\right.\n@space$}}}
\makeatother

\makeatletter
\@addtoreset{equation}{section}
\makeatother
\begin{document}

\newtheorem{thm}{Theorem}
\newtheorem{lem}[thm]{Lemma}
\newdefinition{rmk}{Remark}
\newproof{pf}{Proof}
\newproof{pot}{Proof of Theorem \ref{thm2}}

\begin{frontmatter}



\title{Decay properties of solutions 
to the Cauchy problem 
for the scalar conservation law 
with nonlinearly degenerate viscosity}

\author[rvt]{Natsumi Yoshida}
\ead{14v00067@gst.ritsumei.ac.jp}

\address{BKC Research Organization of Social Sciences, 
Ritsumeikan University, Kusatsu, Shiga 525-8577, Japan
/Osaka City University Advanced Mathematical Institute, 
Sumiyoshi, Osaka 558-8585, Japan.}
\address{}

\begin{abstract}
In this paper, we study the decay rate in time to solutions 
of the Cauchy problem for 
the one-dimensional viscous conservation law 
where the far field states are prescribed.
Especially, we deal with the case that 
the flux function which is convex 
and also the viscosity 
is a nonlinearly degenerate one ($p$-Laplacian type viscosity). 
As the corresponding Riemann problem admits a Riemann solution as 
the constant state or the single rarefaction wave, 
it has already been proved by Matsumura-Nishihara 
that the solution to the Cauchy problem 
tends toward the constant state or 
the single rarefaction wave 
as the time goes to infinity. 
We investigate that the decay rate in time 
of the corresponding solutions. 
Furthermore, we also investigate that the decay rate in time 
of the solution for the higher order derivative. 
These are the first result concerning the asymptotic decay 
of the solutions 
to the Cauchy problem of the scalar conservation law 
with nonlinear viscosity. 
The proof is given by 
$L^1$, $L^2$-energy and time-weighted $L^q$-energy methods. 
\end{abstract}

\end{frontmatter}
%\newtheorem{thm}{Theorem}
%\newtheorem{lem}[thm]{Lemma}
%\newdefinition{rmk}{Remark}
%\newproof{pf}{Proof}
%\newproof{pot}{Proof of Theorem \ref{thm2}}
%

% The thanks line in the title should be filled in if there is
% any support acknowledgement for the overall work to be included
% This \thanks is also used for the received by date info, but
% authors are not expected to provide this.

%This work is 
%supported in part by Grant-in-Aid for Scientific Research (B) 23340036, Japan}}
%         \ \\
% Department of Pure and Applied Mathematics,\\
% Graduate School of Information Science and Technology,\\
% Osaka University, Japan}

%%\begin{AMS}
%%35K55, 35B40, 35L65
%%\end{AMS}

\pagestyle{myheadings}
\thispagestyle{plain}
\markboth{N. YOSHIDA}{NATSUMI YOSHIDA}

\section{Introduction and main theorems}
In this paper, 
we shall consider the asymptotic behavior of 
solutions 
for the one-dimensional scalar conservation law 
with a nonlinearly degenerate viscosity 
($p$-Laplacian type viscosity with $p>1$) 
\begin{eqnarray}
 \left\{\begin{array}{ll}
  \partial_tu +\partial_x \bigl(f(u) \bigr)
  = \mu \, 
    \partial_x \left( \, 
    \left| \, \partial_xu \, \right|^{p-1} \partial_xu \, 
    \right)
  \qquad &(t>0, x\in \mathbb{R}), \\[5pt]
  u(0,x) = u_0(x) \qquad &( x \in \mathbb{R} ),\\[5pt]
  \displaystyle{\lim_{x\to \pm \infty}} u(t,x) =u_{\pm}  
  \qquad &\bigl( t \ge 0 \bigr).   
 \end{array}
 \right.\,
\end{eqnarray}
%for $p>1$. 
Here, $u=u(t,x)$ denotes the unknown function of $t>0$ and $x\in \mathbb{R}$, 
the so-called conserved quantity, 
$f=f(u)$ is the flux function depending only on $u$, 
$\mu$ is the viscosity coefficient, 
$u_0$ is the given initial data, 
and constants $u_{\pm } \in \mathbb{R}$ 
are the prescribed far field states. 
We suppose the given flux $f=f(u)$ is a $C^3$-function % ・ｽ・ｽ ・ｽ・ｽ
satisfying $f(0)=f'(0)=0$, 
$\mu$ is a positive constant 
and far field states $u_{\pm }$ satisfy $u_-<u_+$ 
without loss of generality. 

At first, we shall 
motivate the physical meaning to 
the nonlinearly degenerate viscosity
and review the related models 
conscerning with the Cauchy problem (1.1). 
It is known that if $p=1$ and $f(u)=\frac{1}{2} u^2$, 
the equation in our problem (1.1) becomes 
the viscous Burgers equation: 
\begin{equation*}
\partial_tu +u \, \partial_x u
  = \mu \, 
    \partial_x^2 u. 
\end{equation*}
In particular, the viscosity term $\mu \, \partial_x^2 u$ 
stands for Newtonian fluid. 
The Newtonian fluid is what satisfies the relation between 
the strain rate $\partial_{x_{j}} u_{i} + \partial_{x_{i}} u_{j}$ 
($\partial_x u$, for one-dimensional case) is linear, that is, 
$$
\tau = \mu \left( \, \partial_{x_{j}} u_{i} + \partial_{x_{i}} u_{j} \, \right) 
\quad {\rm{or}} \quad \tau = \mu \, \partial_x u. 
$$
On the other hand, if a fluid satisfies the relation between 
the strain rate and the stress is nonlinear (for example, 
polymers, viscoelastic or viscoplastic flow), 
the fluid is non-Newtonian fluid, 
such as, blood, honey, butter, whipped cream, suspension, and so on. 
The typical nonlinearity in the non-Newtonian fluid 
is the power-law fluid (cf. \cite{ma-pr-st}), that is, 
$$
\tau = \mu \left( \, \partial_{x_{j}} u_{i} + \partial_{x_{i}} u_{j} \, \right)^p 
\quad {\rm{or}} \quad \tau = \mu \left( \, \partial_x u \, \right)^p. 
$$
Lady{\v{z}}enskaja \cite{lad1} %, \cite{lad2}
has proposed a new mathematical model for 
the imcompressible Navier-Stokes equation 
with the power-law type nonlinear viscosity (see also \cite{du-gu}). 
The Lady{\v{z}}enskaja equation is the following: 
$$
\partial_tu_{i} + u_{j} \, \partial_{x_{j}} u_{i} 
= - \partial_{x_{i}} p + 
\partial_{x_{j}} \left( \, \Biggl( \, 
\mu_{0} + \mu_{1} \biggl( \, 
\sum_{\substack{i,j \\ }} \left( \, \partial_{x_{i}} u_{j} \, \right)^2
\, \biggr)^{\frac{r}{2}}
\, \Biggr) \partial_{x_{j}} u_{i} \, \right) 
+ f_{i}
$$
where $i=1,2$, or $i=1,2,3$. 
In particular, if $\mu_{0}=0,\mu_{1}>0$ and $r>-1$, 
this model is said to be the Ostwald-de Waele model: 
$$
\partial_tu_{i} + u_{j} \, \partial_{x_{j}} u_{i} 
= - \partial_{x_{i}} p + 
\partial_{x_{j}} \biggl( \, \Bigl( \, 
\mu \, \left| \, D \vec{u} \, \right|^{r}
\, \Bigr) \partial_{x_{j}} u_{i} \, \biggr) 
+ f_{i}
$$
where 
$\left| \, D \vec{u} \, \right|
:=\biggl( \, 
  \displaystyle{\sum_{\substack{i,j \\ }}} 
  \left( \, \partial_{x_{i}} u_{j} \, \right)^2
  \, \biggr)^{\frac{1}{2}}$, 
and $i=1,2$, or $i=1,2,3$. 
In this sence, our viscosity 
$\mu \, 
    \partial_x \left( \, 
    \left| \, \partial_xu \, \right|^{p-1} \partial_xu \, 
    \right)$ should be called the Ostwald-de Waele type viscosity. 

We are interested in the asymptotic behavior 
and its precise estimates in time of 
the global solution to our problem (1.1). 
It can be expected that the large-time behavior 
is closely related to the weak solution (``Riemann solution'') of the corresponding 
Riemann problem 
(cf. \cite{lax}, \cite{smoller})
for the non-viscous hyperbolic part of (1.1): 
\begin{eqnarray}
 \left\{\begin{array} {ll}
 \partial _t u + \partial _x \bigl( f(u) \bigr)=0 
 \qquad &(t>0, x\in \mathbb{R}),\\[5pt]
u(0,x)=u_0 ^{\rm{R}} (x)\qquad &(x \in \mathbb{R}),
 \end{array}
  \right.\,
\end{eqnarray}
where $u_0 ^{\rm{R}}$ is the Riemann data defined by
$$
u_0 ^{\rm{R}} (x)=u_0 ^{\rm{R}} (x\: ;\: u_- ,u_+)
          := 
          \left\{\begin{array} {ll}
          u_-  & \; (x < 0),\\[5pt]
          u_+  & \; (x > 0).
          \end{array}\right.
$$
In fact, for $p=1$ in (1.1), the usual linear viscosity case: 
\begin{eqnarray}
 \left\{\begin{array}{ll}
  \partial_tu +\partial_x \bigl(f(u) \bigr)= \mu \, \partial_x^2 u
  \qquad &(t>0, x\in \mathbb{R}), \\[5pt]
  u(0,x) = u_0(x) \qquad &( x \in \mathbb{R} ),\\[5pt]
  \displaystyle{\lim_{x\to \pm \infty}} u(t,x) =u_{\pm}  
  \qquad &\bigl( t \ge 0 \bigr),   
 \end{array}
 \right.\,
\end{eqnarray}
when the smooth flux function $f$ is genuinely nonlinear 
on the whole space $\mathbb{R}$,
i.e., $f''(u)\ne 0\ (u\in \mathbb{R})$, 
Il'in-Ole{\u\i}nik \cite{ilin-oleinik} showed 
the following: if  $f''(u) > 0\ (u\in \mathbb{R})$, that is, the Riemann solution 
consists of a single rarefaction wave solution, 
the global solution in time of the Cauchy problem 
(1.3) tends toward the rarefaction wave;
if  $f''(u) < 0\ (u\in \mathbb{R})$, that is,
the Riemann solution consists of a single shock wave solution, 
the global solution of the
Cauchy problem (1.3) does the corresponding smooth traveling wave solution 
(``viscous shock wave'') of (1.3) with a spacial shift 
(cf. \cite{ilin-kalashnikov-oleinik}). 
Hattori-Nishihara \cite{hattori-nishihara} also proved that the asymptotic decay rate in time, 
of the solution toward the single rarefaction wave, is  
$(1+t)^{-{\frac{1}{2}}\, \left(1-\frac{1}{p} \right)}$ 
in the $L^p$-norm $\bigl(1\leq p\leq \infty \bigr)$ for large $t>0$ 
(see also \cite{hashimoto-kawashima-ueda}, \cite{hashimoto-matsumura}, \cite{nakamura}). 
More generally, in the case of the flux functions 
which are not uniformly genuinely nonlinear, 
when the Riemann solution consists of a single shock wave 
satisfying Ole{\u\i}nik's shock condition, 
Matsumura-Nishihara \cite{matsu-nishi3} showed the asymptotic stability 
of the corresponding viscous shock wave. 
Moreover, Matsumura-Yoshida \cite{matsumura-yoshida} 
considered the circumstances 
where the Riemann solution generically forms a pattern of
multiple nonlinear waves which consists of 
rarefaction waves %, shock waves 
and waves of contact discontinuity (refer to \cite{liep-rosh}), 
%there had been no results about the asymptotics toward 
%the multiwave pattern. 
and investigated %focused their attention on 
that the case where the flux function $f$ is smooth
and genuinely nonlinear (that is, $f$ is convex function or concave function) 
on the whole $\mathbb{R}$ except 
a finite interval $I := (a,b) \subset \mathbb{R}$, and 
linearly degenerate on $I$, that is, 
\begin{equation}
\left\{
\begin{array}{ll}
  f''(u) >0 & \; \bigl(u \in (-\infty ,a\, ]\cup [\, b,+\infty )\bigr),\\[5pt]
  f''(u) =0 & \; \bigl(u \in (a,b)\bigr).
\end{array}\right.
\end{equation}
Under the conditions (1.4), 
they proved the unique global solution in time 
to the Cauchy problem (1.3) tends uniformly in space toward 
the multiwave pattern of the combination of 
the viscous contact wave 
%``viscous contact wave for Newtonian'' 
and the rarefaction waves 
as the time goes to infinity. 
Yoshida \cite{yoshida} also obtained 
that the precise decay properties 
for the asymptotics toward the multiwave pattern. 
In fact, owing to \cite{yoshida}, 
the decay rate in time is 
$(1+t)^{-{\frac{1}{2}}\, \left(\frac{1}{2}-\frac{1}{p} \right)}$ 
in the $L^p$-norm $\bigl(2\leq p<+\infty \bigr)$ 
and 
$(1+t)^{-{\frac{1}{4}}\, +\epsilon}$ for any $\epsilon >0$ 
in the $L^\infty$-norm 
if the initial perturbation from the corresponding asymptotics 
satisfies $H^1$. 
Furthermore, 
if the perturbation satisfies $H^1\cap L^1$, 
the decay rate in time is 
$(1+t)^{-{\frac{1}{2}}\, \left(1-\frac{1}{p} \right)\, +\epsilon }$ 
for any $\epsilon >0$ 
in the $L^p$-norm $\bigl(1\leq p<+\infty \bigr)$ 
and
$(1+t)^{-{\frac{1}{2}}\, +\epsilon}$ for any $\epsilon >0$
in the $L^\infty$-norm. 

For $p>1$, 
there are few results for the asymptotic behavior 
for the problem (1.1)
(the related problems are studied in 
\cite{gur-mac}, \cite{nagai-mimura1}, \cite{nagai-mimura2} 
and so on). 
In the case where the flux function is genuinely nonlinear 
on the whole space $\mathbb{R}$, 
Matsumura-Nishihara \cite{matsu-nishi3} 
proved that if the far field states satisfy $u_- = u_+=:\tilde{u}$, 
then the solution tends toward the constant state $\tilde{u}$, 
and if the far field states $u_- < u_+$, 
then the solution tends toward a single rarefaction wave. 
%by using the $L^2$ and $L^p$-energy estimates. 
In the case where the flux function satisfies (1.4), 
Yoshida \cite{yoshida'} recently showed 
that the asymptotics 
which tends toward the multiwave pattern 
of the combination of 
the viscous contact wave 
constructed by the Barenblatt-Kompanceec-Zel'dovi{\v{c}} solution 
(see also \cite{carillo-toscani}, \cite{huang-pan-wang}, \cite{kamin}) 
of the porous medium equation, 
and the rarefaction waves. 
However, the decay rate of any asymptotics of the problem (1.1) 
has been not known.  

The aim of the present paper is to obtain 
the precise time-decay estimates for the asympotics of the 
previous study in \cite{matsu-nishi2}. 

\medskip

\noindent
{\bf Stability Theorem 1.1} (Matsumura-Nishihara \cite{matsu-nishi2}){\bf .}\quad{\it
Let the flux function $f \in C^3(\mathbb{R})$ 
satisfy $f(0)=f'(0)=0$ and $f''(u) > 0\ (u\in \mathbb{R})$,
and the far field states $u_- = u_+=\tilde{u}$. 
Assume that the initial data satisfies %<== Assume (that) ... satisfies ...
$u_0-\tilde{u} \in L^2$ and 
$\partial _x u_0 \in L^{p+1}$. 
Then the Cauchy problem {\rm(1.1)} with $p>1$ has a 
unique global weak 
solution in time $u=u(t,x)$ 
satisfying 
\begin{eqnarray*}
\left\{\begin{array}{ll}
u-\tilde{u} \in C^0\bigl( \, [\, 0,\infty)\, ;L^2 \bigr)
                \cap L^{\infty}\bigl( \, \mathbb{R}^{+} \, ;L^{2} \bigr),\\[5pt]
\partial _x u \in L^{\infty} \bigl( \, \mathbb{R}^{+} \, ;L^{p+1} \bigr)
                  \cap L^{p+1}\bigl(\, {\mathbb{R}^{+}_{t}} \times {\mathbb{R}}_{x} \bigr)
                  \cap L^{p+2}\bigl(\, {\mathbb{R}^{+}_{t}} \times {\mathbb{R}}_{x} \bigr),\\[5pt]
%\partial _t u 
%\in L^{\infty} \bigl( \, \mathbb{R}^{+} \, ;L^{p+1} \bigr)
%,\\[5pt]
\partial_x \left( \, \left| \, \partial_x u \, \right|^{p-1} \partial_xu \, \right)
\in L^{2}\bigl(\, {\mathbb{R}^{+}_{t}} \times {\mathbb{R}}_{x} \bigr) % ・ｽ・ｽ ・ｽ・ｽ・ｽ・ｽ
\end{array} 
\right.\,
\end{eqnarray*}
and the asymptotic behavior 
$$
\lim _{t \to \infty}\sup_{x\in \mathbb{R}}
|\,u(t,x)-\tilde{u}\,| = 0. 
$$
Moreover, 
$\partial _x^2 u 
\in L^{2}\bigl(\, {\mathbb{R}^{+}_{t}} \times {\mathbb{R}}_{x} \bigr)$ 
when $1<p \le \frac{3}{2}$ provided that 
$\partial _x u_0 \in L^{3-p}$. 
}

\medskip

\noindent
{\bf Stability Theorem 1.2} (Matsumura-Nishihara \cite{matsu-nishi2}){\bf .}\quad{\it
Let the flux function $f \in C^3(\mathbb{R})$ 
satisfy $f(0)=f'(0)=0$ and $f''(u) > 0\ (u\in \mathbb{R})$,
and the far field states $u_- < u_+$. 
Assume that the initial data satisfies %<== Assume (that) ... satisfies ...
$u_0-u_0 ^{\rm{R}} \in L^2$ and 
$\partial _x u_0 \in L^{p+1}$. 
Then the Cauchy problem {\rm(1.1)} with $p>1$ has a 
unique global weak 
solution in time $u=u(t,x)$ 
satisfying 
\begin{eqnarray*}
\left\{\begin{array}{ll}
u-u_0 ^{\rm{R}} \in C^0\bigl( \, [\, 0,\infty)\, ;L^2 \bigr)
                \cap L^{\infty}\bigl( \, \mathbb{R}^{+} \, ;L^{2} \bigr),\\[5pt]
\partial _x u \in L^{\infty} \bigl( \, \mathbb{R}^{+} \, ;L^{p+1} \bigr)
                  \cap L^{p+1}\bigl(\, {\mathbb{R}^{+}_{t}} \times {\mathbb{R}}_{x} \bigr)
                  \cap L^{p+2}\bigl(\, {\mathbb{R}^{+}_{t}} \times {\mathbb{R}}_{x} \bigr),\\[5pt]
%\partial _t u 
%\in L^{\infty} \bigl( \, \mathbb{R}^{+} \, ;L^{p+1} \bigr)
%,\\[5pt]
\partial_x \left( \, \left| \, \partial_x u \, \right|^{p-1} \partial_xu \, \right)
\in L^{2}\bigl(\, {\mathbb{R}^{+}_{t}} \times {\mathbb{R}}_{x} \bigr) % ・ｽ・ｽ ・ｽ・ｽ・ｽ・ｽ
\end{array} 
\right.\,
\end{eqnarray*}
and the asymptotic behavior 
$$
\lim _{t \to \infty}\sup_{x\in \mathbb{R}} \, 
\biggl|\,u(t,x)-u^r\left( \, \frac{x}{t}\: ;\: u_-,u_+\right)\,\biggr| = 0, 
$$
where, the rarefaction wave $u^r$ which connects the 
far field states $u_-$ and $u_+$ is explicitly given by 
\begin{equation*}
u^r=u^r \left( \frac{x}{t}\: ;\:  u_- , u_+ \right)
:= \left\{
\begin{array}{ll}
  u_-  & \; \bigl(\, x \leq \lambda(u_-)\,t \, \bigr),\\[7pt]
  \displaystyle{ (\lambda)^{-1}\left( \frac{x}{t}\right) } 
  & \; \bigl(\, \lambda(u_-)\,t \leq x \leq \lambda(u_+)\,t\,  \bigr),\\[7pt]
   u_+ & \; \bigl(\, x \geq \lambda(u_+)\,t \, \bigr),
\end{array}
\right.
\end{equation*} 
where $\lambda(u):=f'(u)$. 
Moreover, 
$\partial _x^2 u 
\in L^{2}\bigl(\, {\mathbb{R}^{+}_{t}} \times {\mathbb{R}}_{x} \bigr)$ 
when $1<p \le \frac{3}{2}$ provided that 
$\partial _x u_0 \in L^{3-p}$. 
}

Now we are ready to state our main results. 

\medskip

\noindent
{\bf Theorem 1.1} (Main Theorem I){\bf .}\quad{\it
Under the same assumptions in Stability Theorem 1.1, 
%for the flux function $f$, the far field states $u_{\pm }\in \mathbb{R}$ 
%and the initial data $u_0$, 
the unique global solution in time $u$ 
of the Cauchy problem {\rm(1.1)} %$u=u(t,x)$ 
satisfying 
\begin{eqnarray*}
\left\{\begin{array}{ll}
u-\tilde{u} \in C^0\bigl( \, [\, 0,\infty)\, ;L^2 \bigr)
                \cap L^{\infty}\bigl( \, \mathbb{R}^{+} \, ;L^{2} \bigr),\\[5pt]
\partial _x u \in L^{\infty} \bigl( \, \mathbb{R}^{+} \, ;L^{p+1} \bigr)
                  \cap L^{p+1}\bigl(\, {\mathbb{R}^{+}_{t}} \times {\mathbb{R}}_{x} \bigr)
                  \cap L^{p+2}\bigl(\, {\mathbb{R}^{+}_{t}} \times {\mathbb{R}}_{x} \bigr),\\[5pt]
%\partial _t u 
%\in L^{\infty} \bigl( \, \mathbb{R}^{+} \, ;L^{p+1} \bigr)
%,\\[5pt]
\partial_x \left( \, \left| \, \partial_x u \, \right|^{p-1} \partial_xu \, \right)
\in L^{2}\bigl(\, {\mathbb{R}^{+}_{t}} \times {\mathbb{R}}_{x} \bigr) % ・ｽ・ｽ ・ｽ・ｽ・ｽ・ｽ
\end{array} 
\right.\,
\end{eqnarray*}
satisfies the following time-decay estimates 
\begin{eqnarray*}
\left\{\begin{array} {ll}
\Arrowvert \, 
u(t) - \tilde{u} 
\, \Arrowvert _{L^q }
\leq C( \, p, \, q, \, u_0 \, ) \, 
     (1+t)^{-\frac{1}{3p+1}\left(1-\frac{2}{q}\right)},\\[5pt]
\Arrowvert \,
u(t) - \tilde{u} 
\, \Arrowvert _{L^{\infty} }
\leq C( \, \epsilon, p, \, q, u_0, \, \partial _x u_0 \, ) \, 
     (1+t)^{-\frac{1}{3p+1}+\epsilon} % ・ｽ・ｽ ・ｽ・ｽ・ｽ・ｽ・ｽ・ｽ・ｽ・ｽ・ｽ・ｽ・ｽ・ｽ
\end{array}
  \right.\,
\end{eqnarray*}
for %$2 \le  p< \infty $
$q \in [\, 2, \infty) $ and any $\epsilon>0$. 
}

\medskip

\noindent
{\bf Theorem 1.2} (Main Theorem I\hspace{-.1em}I){\bf .}\quad{\it
Under the same assumptions in Theorem 1.1, 
%for the flux function $f$, the far field states $u_{\pm }\in \mathbb{R}$ 
%and the initial data $u_0$, 
if the initial data further satisfies $u_0-\tilde{u} \in L^1$, %L^1(\mathbb{R})
then it holds that the unique global solution in time $u$ 
of the Cauchy problem {\rm(1.1)} %$u=u(t,x)$ 
satisfies the following time-decay estimates 
\begin{eqnarray*}
\left\{\begin{array} {ll}
\Arrowvert \,
u(t) - \tilde{u}\,
\Arrowvert _{L^q }
\leq C( \, p, \, q, \, u_0 \, ) \, 
     (1+t)^{-\frac{1}{2p}\left(1-\frac{1}{q}\right)},\\[5pt]
\Arrowvert \,
u(t) - \tilde{u}\,
\Arrowvert _{L^{\infty} }
\leq C( \, \epsilon, \, p, \, q, \, u_0, \, \partial _x u_0 \, ) \, 
     (1+t)^{-\frac{1}{2p}+\epsilon} % ・ｽ・ｽ ・ｽ・ｽ・ｽ・ｽ
\end{array}
  \right.\,
\end{eqnarray*}
for %$1 \le  p< \infty $ 
$q \in [\, 1, \infty) $ and any $\epsilon>0$. 
Furthermore, the solution satisfies 
the following time-decay estimate for the higher order derivative
\begin{eqnarray*}
\bigl|\bigl|\,
\partial _x u(t) \,
\bigr|\bigr|_{L^{p+1} }
\leq C( \, \epsilon , \, p, \, q, \, u_0, \, \partial _x u_0 \, ) \, 
     (1+t)^{-\frac{3}{2(p+1)(3p-2)} + {\epsilon} } % ・ｽ・ｽ ・ｽ・ｽ・ｽ・ｽ・ｽ・ｽ・ｽ・ｽ・ｽ・ｽ・ｽ・ｽ
\end{eqnarray*}
for any $\epsilon>0$. 
}

\medskip

\noindent
{\bf Theorem 1.3} (Main Theorem I\hspace{-.1em}I\hspace{-.1em}I){\bf .}\quad{\it
Under the same assumptions in Theorem 1.2, 
%for the flux function $f$, the far field states $u_{\pm }\in \mathbb{R}$ 
%and the initial data $u_0$, 
if the initial data further satisfies 
$\partial _x u_0 \in L^{r+1} \, (r>p)$, 
then it holds that the unique global solution in time $u$ 
of the Cauchy problem {\rm(1.1)} %$u=u(t,x)$ 
satisfies the following time-decay estimate 
for the higher order derivative
\begin{eqnarray*}
\bigl|\bigl|\,
\partial _x u(t) \,
 \bigr|\bigr|_{L^{r+1} }
\leq C( \, \epsilon , \, p, \, r, \, u_0, \, \partial _x u_0 \, ) \, 
     (1+t)^{-\frac{p+2r}{2p(3p-2)(r+1)} + {\epsilon} } % ・ｽ・ｽ ・ｽ・ｽ・ｽ・ｽ・ｽ・ｽ・ｽ・ｽ・ｽ・ｽ・ｽ・ｽ
\end{eqnarray*}
%\begin{align*}
%\begin{aligned}
%&\bigl|\bigl|\,
% \partial _x u(t) \,
% \bigr|\bigr|_{L^{r+1} } \\
%& \leq \left\{\begin{array} {ll} 
%       C( \, p, \, r, \, \phi_0, \, \partial _x \phi_0 \, ) \, 
%       (1+t)^{-\frac{p+2r}{2p(3p-2)(r+1)}} \, \: \; \quad \qquad 
%       \left( \, 
%       r > \displaystyle{\frac{6p^3+7p^2-2p}{6p^2-2p-2}} 
%       \, \right),\\[10pt] 
%       C( \, p, \, r, \, u_0, \, \partial _x u_0 \, ) \, 
%       (1+t)^{-\frac{6p(r-p)+4(r-p)+3}{2(3p+1)(3p-2)(r+1)}} \: 
%       \left( \, 
%       p < r \le \displaystyle{\frac{6p^3+7p^2-2p}{6p^2-2p-2}} 
%       \, \right).
%       \end{array}
%       \right.\,
%\end{aligned}
%\end{align*}
for any $\epsilon>0$. 
}

\medskip

\noindent
{\bf Theorem 1.4} (Main Theorem I\hspace{-.1em}V){\bf .}\quad{\it
Under the same assumptions in Stability Theorem 1.2, 
%for the flux function $f$, the far field states $u_{\pm }\in \mathbb{R}$ 
%and the initial data $u_0$, 
the unique global solution in time $u$ 
of the Cauchy problem {\rm(1.1)} %$u=u(t,x)$ 
satisfying 
\begin{eqnarray*}
\left\{\begin{array}{ll}
u-u_0 ^{\rm{R}} \in C^0\bigl( \, [\, 0,\infty)\, ;L^2 \bigr)
                \cap L^{\infty}\bigl( \, \mathbb{R}^{+} \, ;L^{2} \bigr),\\[5pt]
\partial _x u \in L^{\infty} \bigl( \, \mathbb{R}^{+} \, ;L^{p+1} \bigr)
                  \cap L^{p+1}\bigl(\, {\mathbb{R}^{+}_{t}} \times {\mathbb{R}}_{x} \bigr)
                  \cap L^{p+2}\bigl(\, {\mathbb{R}^{+}_{t}} \times {\mathbb{R}}_{x} \bigr),\\[5pt]
%\partial _t u 
%\in L^{\infty} \bigl( \, \mathbb{R}^{+} \, ;L^{p+1} \bigr)
%,\\[5pt]
\partial_x \left( \, \left| \, \partial_x u \, \right|^{p-1} \partial_xu \, \right)
\in L^{2}\bigl(\, {\mathbb{R}^{+}_{t}} \times {\mathbb{R}}_{x} \bigr) % ・ｽ・ｽ ・ｽ・ｽ・ｽ・ｽ
\end{array} 
\right.\,
\end{eqnarray*}
satisfies the following time-decay estimates 
\begin{eqnarray*}
\left\{\begin{array} {ll}
\left|\left|\,
u(t) - 
\displaystyle{u^r\left( \, \frac{\cdot }{t}\: ;\: u_-,u_+\right)}\,
\right|\right|_{L^q }
\leq C( \, p, \, q, \, u_0 \, ) \, 
     (1+t)^{-\frac{1}{3p+1}\left(1-\frac{2}{q}\right)},\\[7pt]
\left|\left|\,
u(t) - 
\displaystyle{u^r\left( \, \frac{\cdot }{t}\: ;\: u_-,u_+\right)}
\,\right|\right|_{L^{\infty} }
\leq C( \, \epsilon, \, p, \, q, \, u_0, \, \partial _x u_0 \, ) \, 
     (1+t)^{-\frac{1}{3p+1}+\epsilon} % ・ｽ・ｽ ・ｽ・ｽ・ｽ・ｽ・ｽ・ｽ・ｽ・ｽ・ｽ・ｽ・ｽ・ｽ
\end{array}
  \right.\,
\end{eqnarray*}
for %$2 \le  p< \infty $
$q \in [\, 2, \infty) $ and any $\epsilon>0$. 
}

\medskip

\noindent
{\bf Theorem 1.5} (Main Theorem V){\bf .}\quad{\it
Under the same assumptions in Theorem 1.4, 
%for the flux function $f$, the far field states $u_{\pm }\in \mathbb{R}$ 
%and the initial data $u_0$, 
if the initial data further satisfies $u_0-u_0 ^{\rm{R}} \in L^1$, %L^1(\mathbb{R})
then it holds that the unique global solution in time $u$ 
of the Cauchy problem {\rm(1.1)} %$u=u(t,x)$ 
satisfies the following time-decay estimates 
\begin{eqnarray*}
\left\{\begin{array} {ll}
\left|\left| \,
u(t) - 
\displaystyle{u^r\left( \, \frac{\cdot }{t}\: ;\: u_-,u_+\right)}
\, \right|\right|_{L^q }
\leq C( \, p, \, q, \, u_0 \, ) \, 
     (1+t)^{-\frac{1}{2p}\left(1-\frac{1}{q}\right)},\\[7pt]
\left|\left|\,
u(t) - 
\displaystyle{u^r\left( \, \frac{\cdot }{t}\: ;\: u_-,u_+\right)}
\,\right|\right|_{L^{\infty} }
\leq C( \, \epsilon, \, p, \, q, \, u_0, \, \partial _x u_0 \, ) \, 
     (1+t)^{-\frac{1}{2p}+\epsilon} % ・ｽ・ｽ ・ｽ・ｽ・ｽ・ｽ・ｽ・ｽ・ｽ・ｽ・ｽ・ｽ・ｽ・ｽ
\end{array}
  \right.\,
\end{eqnarray*}
for %$1 \le  p< \infty $ 
$q \in [\, 1, \infty) $ and any $\epsilon>0$. 
Furthermore, the solution satisfies 
the following time-decay estimates for the higher order derivative
\begin{align*}
\begin{aligned}
&\bigl|\bigl|\,
 \partial _x u(t) \,
 \bigr|\bigr|_{L^{p+1} }, \, \; \; 
 \left|\left|\,
 \partial _x u(t) - 
 \displaystyle{\partial _x u^r\left( \, \frac{\cdot }{t}\: ;\: u_-,u_+\right)}
 \,\right|\right|_{L^{p+1} } \\[5pt]
& \leq \left\{\begin{array} {ll} 
      C( \, p, \, u_0, \, \partial _x u_0 \, ) \, 
      (1+t)^{-\frac{p}{p+1}} \, \: \: \; \; \; \quad \qquad 
      \left( \, 1 < p < 
      \displaystyle{\frac{2 + \sqrt{22} }{6} }
      \, \right),\\[15pt] 
      C( \, \epsilon, \, p, \, u_0, \, \partial _x u_0 \, ) \, 
      (1+t)^{-\frac{3}{2(p+1)(3p-2)} + \epsilon} \, \, \: \; \; \; 
      \left( \,  
      \displaystyle{\frac{2 + \sqrt{22} }{6} }
      \leq p \, \right) % ・ｽ・ｽ ・ｽ・ｽ・ｽ・ｽ・ｽ・ｽ・ｽ・ｽ・ｽ・ｽ・ｽ・ｽ
      \end{array}
      \right.\,
\end{aligned}
\end{align*}
for any $\epsilon>0$. 
}

\medskip

\noindent
{\bf Theorem 1.6} (Main Theorem V\hspace{-.1em}I ){\bf .}\quad{\it
Under the same assumptions in Theorem 1.5, 
%for the flux function $f$, the far field states $u_{\pm }\in \mathbb{R}$ 
%and the initial data $u_0$, 
if the initial data further satisfies 
$\partial _x u_0 \in L^{r+1} \, (r>p)$, 
then it holds that the unique global solution in time $u$ 
of the Cauchy problem {\rm(1.1)} %$u=u(t,x)$ 
satisfies the following time-decay estimates 
for the higher order derivative
\begin{align*}
\begin{aligned}
&\bigl|\bigl|\,
 \partial _x u(t) \,
 \bigr|\bigr|_{L^{r+1} }, \, \; \; 
 \left|\left|\,
 \partial _x u(t) - 
 \displaystyle{\partial _x u^r\left( \, \frac{\cdot }{t}\: ;\: u_-,u_+\right)}
 \,\right|\right|_{L^{r+1} } \\[5pt]
& \leq \left\{\begin{array} {ll} 
       C(\, p, \, r, \, u_0, \, \partial _x u_0 \, ) \, 
       (1+t)^{-\frac{2pr+p^2+r}{(3p+1)(r+1)}} \\[15pt] 
       \, \, \, \: \; \; \; \quad \quad \quad \quad \qquad 
       \left( \, 
       1 < p < 
       \displaystyle{\frac{2 + \sqrt{22} }{6} }, \; 
       r > p > \displaystyle{\frac{18p^3-17p^2-16p-3}{2(2p+1)}}
       \, \right),\\[20pt] % ・ｽ・ｽ ・ｽ・ｽ・ｽ・ｽ・ｽ・ｽ・ｽ・ｽ・ｽ・ｽ・ｽ・ｽ
       C( \, \epsilon, \, p, \, r, \, u_0, \, \partial _x u_0 \, ) \, 
       (1+t)^{-\frac{p+2r}{2p(3p-2)(r+1)}+ \epsilon} %\\[15pt] % ・ｽ・ｽ ・ｽ・ｽ・ｽ・ｽ・ｽ・ｽ・ｽ・ｽ・ｽ・ｽ・ｽ・ｽ
%       \, \, \: \; \; \; \quad \quad \quad \quad \quad \qquad 
%       \, \, \; \quad \quad \quad \quad \qquad \qquad \qquad \qquad % ・ｽ・ｽ ・ｽ・ｽ・ｽ・ｽ・ｽ・ｽ・ｽ・ｽ・ｽ・ｽ・ｽ・ｽ
\, \, \: \: \; \; \; \; \; \quad \qquad 
       \left( \, 
       \displaystyle{\frac{2 + \sqrt{22} }{6} }
       \leq p 
       %\le \displaystyle{\frac{2 + \sqrt{7}}{3}}, \; % ・ｽ・ｽ ・ｽ・ｽ・ｽ・ｽ・ｽ・ｽ・ｽ・ｽ・ｽ・ｽ・ｽ・ｽ
%       r > p \ge \displaystyle{\frac{3p^3-2p^2-p}{3p^2-4p-1}}
       \, \right)%,\\[15pt] % ・ｽ・ｽ ・ｽ・ｽ・ｽ・ｽ・ｽ・ｽ・ｽ・ｽ・ｽ・ｽ・ｽ・ｽ
%       C( \, \epsilon, \, p, \, r, \, u_0, \, \partial _x u_0 \, ) \, % ・ｽ・ｽ ・ｽ・ｽ・ｽ・ｽ・ｽ・ｽ・ｽ・ｽ・ｽ・ｽ・ｽ・ｽ
%       (1+t)^{-\frac{6p(r-p)+7p+2r+3}{2(3p+1)(3p-2)(r+1)}+ \epsilon} \\[15pt] 
%       \, \, \, \: \: \; \; \quad \quad \qquad \qquad \qquad \qquad \qquad 
%       \left( \, 
%       p > \displaystyle{\frac{2 + \sqrt{7}}{3}}, \; 
%       p < r \le \displaystyle{\frac{3p^3-2p^2-p}{3p^2-4p-1}}
%       \, \right), % ・ｽ・ｽ ・ｽ・ｽ・ｽ・ｽ・ｽ・ｽ・ｽ・ｽ・ｽ・ｽ・ｽ・ｽ
       \end{array}
       \right.\,
\end{aligned}
\end{align*}
for any $\epsilon>0$. 
}

\medskip

This paper is organized as follows. 
In Section 2, we shall prepare the basic properties of 
the rarefaction wave. 
In Section 3, 
we reformulate the problem in terms of the deviation from 
the asymptotic state 
(similarly in \cite{matsumura-yoshida}, \cite{yoshida}, \cite{yoshida'}), 
that is, 
the single rarefaction wave. 
Following the arguments in \cite{matsu-nishi2}, 
we also prepare some uniform boundedness and energy estimates 
of the deviation 
as the solution to the reformulated problem. 
In order to obtain the time-decay estimates (Theorem 1.4 and Theorem 1.5), 
in Section 4 and Section 5, %<=================・ｽ・ｽ・ｽ・ｽ・ｽ・ｽ・ｽ・ｽ・ｽ・ｽ・ｽ・ｽ・ｽ・ｽ・ｽ・ｽ・ｽ・ｽ・ｽ・ｽ・ｽ\・ｽ・ｽ
we establish the uniform energy estimates in time %<=======・ｽ・ｽ・ｽ・ｽ・ｽ・ｽ・ｽ・ｽ・ｽ・ｽ・ｽ・ｽ
by using a very technical time-weighted energy method. %<===・ｽ・ｽ・ｽ・ｽ・ｽ・ｽ・ｽ・ｽ・ｽ・ｽ・ｽ・ｽ
%and careful estimates of the interactions %<===============・ｽ・ｽ・ｽ・ｽ・ｽ・ｽ・ｽ・ｽ・ｽ・ｽ・ｽ・ｽ
%between the nonlinear waves. 
%<============================・ｽ・ｽ・ｽ・ｽ・ｽ・ｽ・ｽ・ｽ・ｽ・ｽ・ｽ・ｽ
In Section 6, %<==================================・ｽ・ｽ・ｽ・ｽ・ｽ・ｽ・ｽ・ｽ・ｽ・ｽ・ｽ・ｽ
we prove the time-decay $L^{r+1}$-estimate %<========================・ｽ・ｽ・ｽ・ｽ・ｽ・ｽ・ｽ・ｽ・ｽ・ｽ・ｽ・ｽ
for the higher order derivative, Theorem 1.6. %<=================・ｽ・ｽ・ｽ・ｽ・ｽ・ｽ・ｽ・ｽ・ｽ・ｽ・ｽ・ｽ・ｽ・ｽ・ｽ・ｽ・ｽ・ｽ・ｽ・ｽ・ｽ\・ｽ・ｽ
We shall finally discuss the time-decay rates in our main theorems 
comparing with those for a Cauchy problem 
of the symplest $p$-Laplacian evolution equation without convective term in Section 7.

\smallskip

\smallskip

{\bf Some Notation.}\quad 
We denote by $C$ generic positive constants unless 
they need to be distinguished. 
In particular, use $C(\alpha, \beta, \cdots )$ 
or $C_{\alpha, \beta, \cdots }$ 
when we emphasize the dependency on $\alpha, \beta, \cdots $. 
Use $\mathbb{R}^{+}$ as 
$
\mathbb{R}^{+}:=(0,\infty),
$
and the symbol ``$\vee $'' as
$$
a\vee b:= \max \{a,b\}. 
$$
We also use the Friedrichs mollifier $\rho_\delta \ast $, 
where, 
$\rho_\delta(x):=\frac{1}{\delta}\rho \left( \frac{x}{\delta}\right)$
with 
\begin{align*}
\begin{aligned}
&\rho \in C^{\infty}_0(\mathbb{R}),\quad 
\rho (x)\geq 0\:  (x \in \mathbb{R}), \\
&\mathrm{supp} \{\rho \} \subset 
\left\{x \in \mathbb{R}\: \left|\:  |\, x \, |\le 1 \right. \right\},\quad  
\int ^{\infty}_{-\infty} \rho (x)\, \mathrm{d}x=1, 
\end{aligned}
\end{align*}
and $\rho_\delta \ast f$ denote the convolution. 
For function spaces, 
$L^p = L^p(\mathbb{R})$ and $H^k = H^k(\mathbb{R})$ 
denote the usual Lebesgue space and 
$k$-th order Sobolev space on the whole space $\mathbb{R}$ 
with norms $||\cdot||_{L^p}$ and $||\cdot||_{H^k}$, respectively. 
We also define 
the bounded $C^{m}$-class $\mathscr{B}^{m}$ as follows 
$$
f\in \mathscr{B}^{m}(\Omega)
\, \Longleftrightarrow  \, 
f\in C^m(\Omega), 
\; 
\sup _{\Omega}\, \sum _{k=0}^{m} \, \bigl| \, D^kf\, \bigr|<\infty 
$$
for $m< \infty$ and 
$$
f\in \mathscr{B}^{\infty }(\Omega)
\, \Longleftrightarrow  \, 
\forall n\in \mathbb{N},\, f\in C^n(\Omega), 
\; 
\sup _{\Omega}\, \sum _{k=0}^{n} \, \bigl| \, D^kf\, \bigr|<\infty 
$$
where $\Omega \subset \mathbb{R}^d$ and 
$D^k$ denote the all of $k$-th order derivatives. 

\section{Preliminaries} 
In this section, 
we shall arrange the two lemmas concerned with % ・ｽ・ｽ ・ｽ・ｽ
the basic properties of 
the rarefaction wave 
for accomplishing the proof of our main theorems. 
Since the rarefaction wave $u^r$ is not smooth enough, 
we need some smooth approximated one 
as in the previous results in 
\cite{hashimoto-matsumura}, \cite{liu-matsumura-nishihara}, \cite{matsu-nishi1}, \cite{matsumura-yoshida}. 
We start with the well-known arguments on $u^r$ 
and the method of constructing its smooth approximation. 
We first consider the rarefaction wave solution $w^r$ %in distribution sense 
to the Riemann problem 
for the non-viscous Burgers equation %(hyperbolic equation)
%with prescribed far field states $w_\pm \in \mathbb{R} \: (w_-<w_+)$: 
\begin{equation}
\label{riemann-burgers}
  \left\{\begin{array}{l}
  \partial _t w + 
  \displaystyle{ \partial _x \left( \, \frac{1}{2} \, w^2 \right) } = 0 
  \, \, \; \; \qquad \quad \qquad ( t > 0,\,x\in \mathbb{R}),\\[7pt]
  w(0,x) = w_0 ^{\rm{R}} ( \, x\: ;\: w_- ,w_+):= \left\{\begin{array}{ll}
                                                  w_+ & \, \: \; \quad (x>0),\\[5pt]
                                                  w_- & \, \: \; \quad (x<0),
                                                  \end{array}
                                                  \right.
  \end{array}
  \right.
\end{equation}
where $w_\pm \in \mathbb{R} \: (w_-<w_+)$ are 
the prescribed far field states. 
The unique global weak solution 
$w=w^r\left( \, \frac{x}{t}\: ;\: w_-,w_+\right)$ 
of (\ref{riemann-burgers}) is explicitly given by 
\begin{equation}
\label{rarefaction-burgers}
w^r \left( \, \frac{x}{t}\: ;\: w_-,w_+\right) := 
  \left\{\begin{array}{ll}
  w_{-} & \bigl(\, x \leq w_{-}t \, \bigr),\\[5pt]
  \displaystyle{ \frac{x}{t} } & \bigl(\, w_{-} t \leq x \leq w_{+} t \, \bigr),\\[5pt]
  w_+ & \bigl(\, x\geq w_{+} t \, \bigr).
  \end{array}\right.
\end{equation} 
Next, under the condition 
$f''(u)>0\ (u\in \mathbb{R})$ and $u_-<u_+$, 
the rarefaction wave solution 
$u=u^r\left( \, \frac{x}{t}\: ;\: u_-,u_+\right)$ 
of the Riemann problem (1.2) 
for hyperbolic conservation law 
is exactly given by 
\begin{equation}
u^r\left( \, \frac{x}{t} \: ; \:  u_-,u_+\right) 
= (\lambda)^{-1}\biggl( w^r\left( \, \frac{x}{t} \: ; \:  \lambda_-,\lambda_+\right)\biggr)
\end{equation}
which is nothing but (1.6), 
where $\lambda_\pm := \lambda(u_\pm) = f'(u_\pm)$. 
%We make a smooth approximation of the rarefaction wave $u^r$ 
%due to the previous works 
%(cf. \cite{hashimoto-matsumura},\cite{liu-matsumura-nishihara},\cite{matsu-nishi1},\cite{matsumura-yoshida}). 
We define a smooth approximation of $w^r(\, \frac{x}{t}\: ;\: w_-,w_+)$ 
by the unique classical solution 
$$
w=w(\, t,x\: ;\: w_-,w_+)\in \mathscr{B}^{\infty }( \, [\, 0,\infty )\times \mathbb{R})
$$
to the Cauchy problem for the following %hyperbolic 
non-viscous Burgers equation
\begin{eqnarray}
\label{smoothappm}
\left\{\begin{array}{l}
 \partial _t w 
 + \displaystyle{ \partial _x \left( \, \frac{1}{2} \, w^2 \right) } =0 
 \, \, \; \; \quad \qquad \qquad \qquad \qquad \qquad  
 (\ t>0,\,x\in \mathbb{R}),\\[7pt]
 w(0,x) 
 = w_0(x) 
 := \displaystyle{ \frac{w_-+w_+}{2} + \frac{w_+-w_-}{2}\tanh x } 
 \qquad \quad \; \:  (x\in \mathbb{R}),
\end{array}
\right.
\end{eqnarray}   
%Then, we define the approximation for 
By using the method of characteristics, 
we get the following formula
\begin{eqnarray}
 \left\{\begin{array} {l}
 w(t,x)=w_0\bigl(x_0(t,x)\bigr)=
 \displaystyle{ \frac{\lambda_-+\lambda_+}{2} } 
+ \displaystyle{ \frac{\lambda_+-\lambda_-}{2}\tanh \bigl( x_0(t,x) \bigr)} ,\\[7pt]
 x=x_0(t,x)+w_0\bigl(x_0(t,x)\bigr)\,t.
 \end{array}
  \right.\,
\end{eqnarray}
We also note the assumption of the flux function $f$ to be % ・ｽ・ｽ ・ｽ・ｽ
$\lambda'(u)\left( =\frac{\mathrm{d}^2f}{\mathrm{d}u^2}(u)\right)>0$. 
%%%%%%%%%%%%%%%%%%%%%%%%%%%%%%%%%%%%%%%%%%%%%%%%%%%%%%%%%%%%%%%%%%%%%%%%%%%%%%%%%%%・ｽ・ｽ
%%%%%%%%%%%%%%%                     w(t,x)                    %%%%%%%%%%%%%%%%%%%%%・ｽ・ｽ
%%%%%%%%%%%%%%%%%%%%%%%%%%%%%%%%%%%%%%%%%%%%%%%%%%%%%%%%%%%%%%%%%%%%%%%%%%%%%%%%%%%・ｽ・ｽ

Now we summarize the results for the smooth approximation $w(\, t,x\: ;\: w_-,w_+)$ 
in the next lemma. %in Lemma 2.1. 
Since the proof is given by the direct calculation as in \cite{matsu-nishi1}, %of (2.5). 
we omit it. %the details. 

\medskip

\noindent
%     w\bigl(\, t,x\: ;\: f^{\epsilon}(w_-), \, f^{\epsilon}(w_+) \, \bigr)
{\bf Lemma 2.1.}\quad{\it
Assume that the far field states satisfy $w_-<w_+$. 
Then the classical solution $w(t,x)=w(\, t,x\: ;\: w_-,w_+)$
given by {\rm(2.4)} 
satisfies the following properties: 

\noindent
{\rm (1)}\ \ $w_- < w(t,x) < w_+$ and\ \ $\partial_xw(t,x) > 0$  
\quad  $(t>0, x\in \mathbb{R})$.

\smallskip

\noindent
{\rm (2)}\ For any $1\leq q \leq \infty$, there exists a positive 
constant $C_q$ such that
             \begin{eqnarray*}
                 \begin{array}{l}
                    \parallel \partial_x w(t)\parallel_{L^q} \leq 
                    C_q (1+t)^{-1+\frac{1}{q}} 
                    \; \quad \bigl(t\ge 0 \bigr),\\[5pt]
                    \parallel \partial_x^2 w(t) \parallel_{L^q} \leq 
                    C_q (1+t)^{-1} 
                    \, \, \quad \quad \bigl(t\ge 0 \bigr).
                    \end{array}       
              \end{eqnarray*}
              
\smallskip

\noindent
{\rm (3)}\ $\: \displaystyle{\lim_{t\to \infty} 
\sup_{x\in \mathbb{R}}
\left| \,w(t,x)- w^r \left( \frac{x}{t} \right) \, \right| = 0}.$
}

\bigskip
%%%%%%%%%%%%%%%%%%%%%%%%%%%%%%%%%%%%%%%%%%%%%%%%%%%%%%%%%%%%%%%%%%%%%%%%%%%%%%%%%%%・ｽ・ｽ
%%%%%%%%%%%%%%%                     w(t,x)                    %%%%%%%%%%%%%%%%%%%%%・ｽ・ｽ
%%%%%%%%%%%%%%%%%%%%%%%%%%%%%%%%%%%%%%%%%%%%%%%%%%%%%%%%%%%%%%%%%%%%%%%%%%%%%%%%%%%・ｽ・ｽ

\noindent
We define the approximation for 
the rarefaction wave $u^r\left( \, \frac{x}{t}\: ;\: u_-,u_+\right)$ by 
\begin{equation}
U^r(\, t,x\: ; \: u_-,u_+) := (\lambda)^{-1} \bigl( w(\, t,x\: ;\: \lambda_-,\lambda_+)\bigr).
\end{equation}

Then we have the next lemma 
as in the previous works 
(cf. \cite{hashimoto-matsumura}, \cite{liu-matsumura-nishihara}, \cite{matsu-nishi1}, \cite{matsumura-yoshida}). 

\medskip

\noindent
{\bf Lemma 2.2.}\quad{\it% <----- Lemma 2.2・ｽ・ｽ・ｽ・ｽ・ｽ・ｽ・ｽ・ｽ・ｽ・ｽ・ｽ・ｽ・ｽ・ｽ・ｽ・ｽ・ｽ・ｽ・ｽ・ｽ・ｽ・ｽ・ｽ・ｽ・ｽ・ｽ・ｽ・ｽ・ｽ・ｽ・ｽ・ｽ
Assume that the far field states satisfy $u_-<u_+$, 
and the flux fanction $f\in C^3(\mathbb{R})$, $f''(u)>0 \: (u\in [\,u_-,u_+\,])$. 
Then we have the following properties:

\noindent
{\rm (1)}\ $U^r(t,x)$ defined by {\rm (2.6)} is %<----- 2.6・ｽ・ｽ・ｽ・ｽ・ｽ・ｽ・ｽ・ｽ・ｽ・ｽ・ｽ・ｽ・ｽ・ｽ・ｽ・ｽ・ｽ・ｽ・ｽ・ｽ・ｽ・ｽ・ｽ・ｽ・ｽ・ｽ・ｽ・ｽ・ｽ・ｽ・ｽ・ｽ
the unique $C^2$-global solution in space-time 
of the Cauchy problem
%\begin{equation*}
$$
\left\{
\begin{array}{l} 
\partial _t U^r +\partial _x \bigl( f(U^r ) \bigr) = 0 
\, \, \, \, \; \; \; \quad \quad \qquad \qquad \qquad \qquad  
(t>0, x\in \mathbb{R}),\\[7pt]
U^r(0,x) 
= \displaystyle{ (\lambda)^{-1} \left( \, \frac{\lambda_- + \lambda_+}{2} 
+ \frac{\lambda_+ - \lambda_-}{2} \tanh x \, \right) } 
\; \: \: \quad \quad( x\in \mathbb{R}),\\[7pt]
\displaystyle{\lim_{x\to \pm \infty}} U^r(t,x) =u_{\pm} 
\, \, \: \: \; \; \quad \quad \qquad \qquad \qquad \qquad \qquad \qquad 
\bigl(t\ge 0 \bigr).
\end{array}
\right.\,     
$$
%\end{equation*}
{\rm (2)}\ \ $u_- < U^r(t,x) < u_+$ and\ \ $\partial_xU^r(t,x) > 0$  
\quad  $(t>0, x\in \mathbb{R})$.

\smallskip

\noindent
{\rm (3)}\ For any $1\leq q \leq \infty$, there exists a positive 
constant $C_q$ such that
             \begin{eqnarray*}
                 \begin{array}{l}
                    \parallel \partial_x %U^r(t,\cdot \: )
                    U^r(t) 
                    \parallel_{L^q} \leq 
                    C_q(1+t)^{-1+\frac{1}{q}} 
                    \quad \bigl(t\ge 0 \bigr),\\[5pt]
                    \parallel \partial_x^2 U^r(t) \parallel_{L^q} \leq 
                    C_q(1+t)^{-1}
                    \, \, \quad \quad \bigl(t\ge 0 \bigr).
                    \end{array}       
              \end{eqnarray*}
              
\smallskip

\noindent
{\rm (4)}\ $\: \displaystyle{\lim_{t\to \infty} 
\sup_{x\in \mathbb{R}}
\left| \,U^r(t,x)- u^r \left( \frac{x}{t} \right) \, \right| = 0}.$

\smallskip

\noindent
{\rm (5)}\ For any $\epsilon \in (0,1)$, there exists a positive 
constant $C_\epsilon$ such that
$$
\left|\, 
U^r(t,x)-u_+ \, \right|
\leq C_\epsilon(1+t)^{-1+\epsilon}\mathrm{e}^{-\epsilon |x-\lambda_+t|}
\quad \bigl(t\ge 0, x \ge \lambda_+t \bigr).
$$

\smallskip

\noindent
{\rm (6)}\ For any $\epsilon \in (0,1)$, there exists a positive 
constant $C_\epsilon$ such that
$$
\left|\, 
U^r(t,x)-u_- \, \right|
\leq C_\epsilon(1+t)^{-1+\epsilon}\mathrm{e}^{-\epsilon |x-\lambda_-t|}
\quad \bigl(t\ge 0, x \le \lambda_-t\bigr).
$$

\smallskip

\noindent
{\rm (7)}\ For any $\epsilon \in (0,1)$, there exists a positive 
constant $C_\epsilon$ such that
$$
\left| \,U^r(t,x) - u^r\left( \frac{x}{t}\right) \, \right| 
\leq C_\epsilon(1+t)^{-1+\epsilon} 
\qquad \bigl(t \ge 1,  \lambda_-t \le x \le \lambda_+t\bigr).
$$

\smallskip

\noindent
{\rm (8)}\ For any $(\epsilon ,q)\in (0,1)\times [\, 1,\infty \, ]$, 
there exists a positive 
constant $C_{\epsilon,q}$ such that
$$
\left|\left|
 \,U^r(t,\cdot \: ) - u^r\left( \frac{\cdot }{t}\right) \, 
\right|\right|_{L^q}  
\leq C_{\epsilon,q}(1+t)^{-1+\frac{1}{q}+\epsilon} 
\qquad \bigl(t \ge 0 \bigr).
$$
}

\noindent
Because the proofs of (1) to (4) are given in \cite{matsu-nishi1}, 
(5) to (7) are in \cite{matsumura-yoshida} 
and (8) is in \cite{yoshida}, 
we omit the proofs here.

%%%%%%%%%%%%%%%%%%%%%%%%%%%%%%%%%%%%%%%%%%%%%%%%%%%%%%%%%%%%%%%%%%%%%%%%%%%%%%%%%%%%%%%%%%%%
%%%%%%%%%%%%%%%%%%%%%%%%%%%%%%%%%%%%%%%%%%%%%%%%%%%%%%%%%%%%%%%%%%%%%%%%%%%%%%%%%%%%%%%%%

\section{Reformulation of the problem}
In this section, 
we reformulate our Cauchy problem (1.1)
in terms of the deviation from 
the asymptotic state, 
the single rarefaction wave. 
We first should note by Lemma 2.2, 
the asymptotic state $u^r\left( \, \frac{x}{t}\: ;\: u_-,u_+\right)$ 
can be replaced by 
\begin{equation*}
U^r( t,x \: ; \: u_-, u_+). 
\end{equation*}
In fact, 
from Lemma 2.1 (especially (8)), 
it follows that 
for any $\epsilon>0$
\begin{align*}
\begin{aligned}
&\left|\left| \, 
 U^r( t,\cdot \: ; \: u_-, u_+)
 - u^r\left( \, \frac{\cdot}{t}\: ;\: u_-,u_+\right)\, 
 \right|\right|_{L^q}\\
&\le C_{\epsilon,q}(1+t)^{-\left( 1-\frac{1}{q}\right)+\epsilon}
\qquad (t\ge 0\, ;\, 1\le q\le \infty).
\end{aligned}
\end{align*}
Then it is noted that $U^r$ is 
monotonically increasing and approximately satisfies 
the equation of (1.1) as
\begin{equation}
 \partial _tU^r +\partial_x \bigl(f(U^r)\bigr)= 0.
\end{equation}
Now  putting 
\begin{equation}
u(t,x) = U^r(t,x) + \phi(t,x)
\end{equation}
and using (3.1),
we can reformulate the problem (1.1) in terms of 
the deviation $\phi $ from 
$U^r$ as 
\begin{eqnarray}
 \left\{\begin{array}{ll}
  \partial _t\phi + \partial_x \left( f(U^r+\phi) - f(U^r) \right) \\[5pt]
  - \mu \, \partial_x 
    \left( \, 
    \bigl| \, \partial_x U^r + \partial_x \phi \, \bigr|^{p-1} 
    \bigl( \, \partial_x U^r + \partial_x \phi \, \bigr) - 
    \bigl| \, \partial_x U^r \, \bigr|^{p-1} \partial_x U^r  \, 
    \right)\\[2pt]
    \, \: \; \quad \quad \qquad \qquad 
    = \mu \, 
      \partial_x \left( \, 
      \bigl| \, \partial_x U^r \, \bigr|^{p-1} \partial_x U^r  
      \, \right)
  \, \, \, \: \: \quad  (t>0, x\in \mathbb{R}), \\[5pt]
  \phi(0,x) = \phi_0(x) 
  := u_0(x)-U^r(0,x) 
  \qquad \qquad \quad \; \: \: \; \; \: \; \; \,\, \, \, (x\in \mathbb{R}). 
 \end{array}
 \right.\,
\end{eqnarray}
Then we look for the global solution in time 
$$
\phi \in C^0\bigl( \, [\, 0,\infty)\, ;L^2 \bigr)
         \cap L^{\infty}\bigl( \, \mathbb{R}^{+} \, ;L^{2} \bigr) 
$$ with 
$$
\partial _x \phi \in L^{\infty} \bigl( \, \mathbb{R}^{+} \, ;L^{p+1} \bigr)
                     \cap L^{p+1}
                     \bigl(\, {\mathbb{R}^{+}_{t}} 
                     \times {\mathbb{R}}_{x} \bigr). 
$$
Here we note that the assumptions on $u_0$ 
and the fact 
$
\partial_x U^r(0,\cdot \, )\in L^{p+1}
$ 
imply $\phi_0 \in L^2$ with $\partial _x \phi_0 \in L^{p+1}$. 
Then the corresponding our main theorems for $\phi$ 
we should prove are as follows. 

\medskip

\noindent
{\bf Theorem 3.1.}\quad{\it
Assume that the flux function 
$f\in C^3 (\, [\, 0,\infty))$ 
satisfies $f(0)=f'(0)=0$ and $f''(u) > 0\ (u\in \mathbb{R})$,  
the far field states $u_-<0<u_+$, 
and the initial data 
$\phi_0 \in L^2$ and 
$\partial _x u_0 \in L^{p+1}$. 
Then, the unique global solution in time $\phi$ 
of the Cauchy problem {\rm(3.3)} 
satisfying 
\begin{eqnarray*}
\left\{\begin{array}{ll}
\phi \in C^0\bigl( \, [\, 0,\infty)\, ;L^2 \bigr)
         \cap L^{\infty}\bigl( \, \mathbb{R}^{+} \, ;L^{2} \bigr),\\[5pt]
\partial _x \phi \in L^{\infty} \bigl( \, \mathbb{R}^{+} \, ;L^{p+1} \bigr)
                     \cap L^{p+1}\bigl(\, {\mathbb{R}^{+}_{t}} \times {\mathbb{R}}_{x} \bigr),\\[5pt]
\partial_x (U^r + \phi) \in L^{\infty} \bigl( \, \mathbb{R}^{+} \, ;L^{p+1} \bigr)
                            \cap L^{p+1}\bigl(\, {\mathbb{R}^{+}_{t}} \times {\mathbb{R}}_{x} \bigr)
                            \cap L^{p+2}\bigl(\, {\mathbb{R}^{+}_{t}} \times {\mathbb{R}}_{x} \bigr),\\[5pt]
%\partial _t u 
%\in L^{\infty} \bigl( \, \mathbb{R}^{+} \, ;L^{p+1} \bigr)
%,\\[5pt]
\partial_x 
\left( \, \left| \, \partial_x (U^r + \phi) \, \right|^{p-1} 
\partial_x(U^r + \phi) \, \right)
\in L^{2}\bigl(\, {\mathbb{R}^{+}_{t}} \times {\mathbb{R}}_{x} \bigr)
\end{array} 
\right.\,
\end{eqnarray*}
and 
$$
\displaystyle{\lim_{t\to \infty}}\sup_{x\in \mathbb{R}}\, 
\bigl| \, \phi(t,x)\,\bigr| = 0
$$
satisfies the following time-decay estimates 
\begin{eqnarray*}
\left\{\begin{array} {ll}
\left|\left|\,
\phi (t) 
\, \right|\right|_{L^q }
\leq C( \, p, \, q, \, \phi_0 \, ) \, 
     (1+t)^{-\frac{1}{3p+1}\left(1-\frac{2}{q}\right)},\\[7pt]
\left|\left|\,
\phi (t) 
\, \right|\right|_{L^{\infty} }
\leq C( \, \epsilon, \, p, \, q, \, \phi_0, \, \partial _x u_0 \, ) \, 
     (1+t)^{-\frac{1}{3p+1}+\epsilon} % ・ｽ・ｽ ・ｽ・ｽ・ｽ・ｽ・ｽ・ｽ・ｽ・ｽ・ｽ・ｽ・ｽ・ｽ
\end{array}
  \right.\,
\end{eqnarray*}
for %$2 \le  p< \infty $
$q \in [\, 2, \infty) $ and any $\epsilon>0$. 
}

\medskip

\noindent
{\bf Theorem 3.2.}\quad{\it
Under the same assumptions in Theorem 3.1, 
if the initial data further satisfies $\phi_0 \in L^1$, %L^1(\mathbb{R})
then it holds that the unique global solution in time $\phi$ 
of the Cauchy problem {\rm(3.3)} 
satisfies the following time-decay estimates 
\begin{eqnarray*}
\left\{\begin{array} {ll}
\left|\left| \,
\phi (t) 
\, \right|\right|_{L^q }
\leq C( \, p, \, q, \, \phi_0 \, ) \, 
     (1+t)^{-\frac{1}{2p}\left(1-\frac{1}{q}\right)},\\[7pt]
\left|\left|\,
\phi (t) 
\,\right|\right|_{L^{\infty} }
\leq C( \, \epsilon, \, p, \, q, \, \phi_0, \, \partial _x u_0 \, ) \, 
     (1+t)^{-\frac{1}{2p}+\epsilon} % ・ｽ・ｽ ・ｽ・ｽ・ｽ・ｽ・ｽ・ｽ・ｽ・ｽ・ｽ・ｽ・ｽ・ｽ
\end{array}
  \right.\,
\end{eqnarray*}
for %$1 \le  p< \infty $ 
$q \in [\, 1, \infty) $ and any $\epsilon>0$. 
Furthermore, the solution satisfies 
the following time-decay estimates for the higher order derivative
\begin{align*}
\begin{aligned}
&\bigl|\bigl|\,
 \partial _x u(t) \,
 \bigr|\bigr|_{L^{p+1} }, \, \; \; 
 \left|\left|\,
 \partial _x \phi (t) 
 \,\right|\right|_{L^{p+1} } \\[5pt]
& \leq \left\{\begin{array} {ll} 
      C( \, \epsilon, \, p, \, \phi_0, \, \partial _x u_0 \, ) \, 
      (1+t)^{-\frac{p}{p+1}} \, \: \; \qquad 
      \left( \, 1 < p \le 
      \displaystyle{\frac{1}{3} + 
      \sqrt{\frac{11}{18} - \frac{(p+1)(3p-2)}{3}\, \epsilon } }
      \, \right),\\[15pt] 
      C( \, \epsilon, \, p, \, \phi_0, \, \partial _x u_0 \, ) \, 
      (1+t)^{-\frac{3}{2(p+1)(3p-2)} + \epsilon} \, \, \; 
      \left( \,  
      \displaystyle{\frac{1}{3} + 
      \sqrt{\frac{11}{18} - \frac{(p+1)(3p-2)}{3}\, \epsilon } }
      < p \, \right) % ・ｽ・ｽ ・ｽ・ｽ・ｽ・ｽ・ｽ・ｽ・ｽ・ｽ・ｽ・ｽ・ｽ・ｽ
      \end{array}
      \right.\,
\end{aligned}
\end{align*}
for any $0<\epsilon \ll 1$. 
}

\medskip

\noindent
{\bf Theorem 3.3.}\quad{\it
Under the same assumptions in Theorem 3.2, 
if the initial data further satisfies 
$\partial _x u_0 \in L^{r+1} \, (r>p)$, 
then it holds that the unique global solution in time $\phi$ 
of the Cauchy problem {\rm(3.3)} 
satisfies the following time-decay estimates 
for the higher order derivative
\begin{align*}
\begin{aligned}
&\bigl|\bigl|\,
 \partial _x u(t) \,
 \bigr|\bigr|_{L^{r+1} }, \, \; \; 
 \left|\left|\,
 \partial _x \phi(t) 
 \,\right|\right|_{L^{r+1} } \\[5pt]
& \leq \left\{\begin{array} {ll} 
       C( \, \epsilon, \, p, \, r, \, \phi_0, \, \partial _x u_0 \, ) \, 
       (1+t)^{-\frac{2pr+p^2+r}{(3p+1)(r+1)}} \\[10pt] 
       \, \, \: \; \; 
       \left( \, 
       1 < p \le 
       \displaystyle{\frac{1}{3} + 
      \sqrt{\frac{11}{18} - \frac{p(3p-2)(r+1)}{2(r-p+1)}\, \epsilon } }, \; 
       r > p > \displaystyle{\frac{18p^3-17p^2-16p-3}{2(2p+1)}}
       \, \right),\\[15pt] 
       C( \, \epsilon, \, p, \, r, \, \phi_0, \, \partial _x u_0 \, ) \, 
       (1+t)^{-\frac{p+2r}{2p(3p-2)(r+1)}+ \epsilon} \\[15pt] 
       \, \, \, \, \: \; \; \; \quad \quad \quad \qquad \qquad \qquad \qquad \qquad \qquad 
       \left( \, 
       \displaystyle{\frac{1}{3} + 
      \sqrt{\frac{11}{18} - \frac{p(3p-2)(r+1)}{2(r-p+1)}\, \epsilon } }
       < p 
       %\le \displaystyle{\frac{2 + \sqrt{7}}{3}}, \;  % ・ｽ・ｽ ・ｽ・ｽ・ｽ・ｽ・ｽ・ｽ・ｽ・ｽ・ｽ・ｽ・ｽ・ｽ
%       r > p \ge \displaystyle{\frac{3p^3-2p^2-p}{3p^2-4p-1}} % ・ｽ・ｽ ・ｽ・ｽ・ｽ・ｽ・ｽ・ｽ・ｽ・ｽ・ｽ・ｽ・ｽ・ｽ
       \, \right)%,\\[15pt]  % ・ｽ・ｽ ・ｽ・ｽ・ｽ・ｽ・ｽ・ｽ・ｽ・ｽ・ｽ・ｽ・ｽ・ｽ
%       C( \, \epsilon, \, p, \, r, \, \phi_0, \, \partial _x u_0 \, ) \, 
%       (1+t)^{-\frac{6p(r-p)+7p+2r+3}{2(3p+1)(3p-2)(r+1)}+ \epsilon} \\[10pt] 
%       \, \, \, \, \: \: \; \; \; \; \quad \quad \quad \quad \qquad \qquad \qquad \qquad \qquad 
%       \left( \, 
%       p > \displaystyle{\frac{2 + \sqrt{7}}{3}}, \; 
%       p < r \le \displaystyle{\frac{3p^3-2p^2-p}{3p^2-4p-1}}
%       \, \right), % ・ｽ・ｽ ・ｽ・ｽ・ｽ・ｽ・ｽ・ｽ・ｽ・ｽ・ｽ・ｽ・ｽ・ｽ
       \end{array}
       \right.\,
\end{aligned}
\end{align*}
for any $0<\epsilon \ll 1$. 
}

\medskip
In order to accomplish the proofs of 
Theorem 3.1, Theorem 3.2 and Theorem 3.3, 
we will need some estimates about boundedness 
of the perturbation $\phi$ and $u$.  
We shall arrange some lemmas for them. 

By using the maximum principle 
(cf. \cite{ilin-kalashnikov-oleinik}, \cite{ilin-oleinik}), 
we first have the following uniform boundedness 
of the perturbation $\phi$ (and also $u$), that is, 

\medskip

\noindent
{\bf Lemma 3.1} (uniform boundedness){\bf .}\quad {\it
It holds that 
\begin{equation}
\sup_{t\in [\, 0,\infty), x\in \mathbb{R}}|\, \phi(t,x)\,| 
\le \| \,\phi_0\, \|_{L^{\infty}} 
    + 2 \, \bigl(\, | \, u_- \, | + | \, u_+ \, | \, \bigr), 
\end{equation}
\begin{align}
\begin{aligned}
&\sup_{t\in [\, 0,\infty), x\in \mathbb{R}}|\, u(t,x)\,| \\
&\le \| \,\phi_0\, \|_{L^{\infty}} 
     + 2 \, \bigl(\, | \, u_- \, | + | \, u_+ \, | \, \bigr)
     + | \, u_- \, | \vee | \, u_+ \, | 
     =: \widetilde{C}. 
\end{aligned}
\end{align}
}

\medskip

\noindent
Secondly, 
we also have the uniform estimates of $\phi$ as follows
(for the proof of it, 
see \cite{matsu-nishi2} ). 

\medskip

\noindent
{\bf Lemma 3.2} (uniform estimates){\bf .}\quad {\it
The unique global solution in time $\phi$ of the Cauchy problem {\rm (3.3)} 
satisfies the following uniform energy inequalities 
\begin{equation}
\| \, \phi (t) \, \| _{L^2}^2
+\int _0^{\infty}
  \| \, \partial _x\phi (t) \, \|_{L^{p+1}}^{p+1} 
  \, \mathrm{d}t 
\leq C_{p}\bigl(\, \| \, \phi _0 \, \| _{L^2} \, \bigr), 
\end{equation}
\begin{align}
\begin{aligned}
\| \, \partial _x u (t) 
\, \| _{L^{p+1}}^{p+1} 
&+\int _0^{\infty} \int _{-\infty}^{\infty} 
 \bigl| \, \partial _x u \, \bigr|^{2(p-1)} 
 \left( \, \partial _x^2 u \, \right)^2 
 \, \mathrm{d}x \mathrm{d}t \\
&\qquad \quad \leq 
 C_{p}\bigl(\, \| \, \phi _0 \, \| _{L^2}, 
 \| \, \partial _x u_0 \, \| _{L^{p+1}} \, \bigr), 
\end{aligned}
\end{align}
\begin{equation}
\int _0^{\infty}
 \| \, 
 \partial _x u (t) 
 \, \| _{L^{p+2}}^{p+2} 
 \, \mathrm{d}t 
\leq C_{p}\bigl(\, \| \, \phi _0 \, \| _{L^2}, 
 \| \, \partial _x u_0 \, \| _{L^{p+1}} \, \bigr)
\end{equation}
for $t \in [\, 0,\infty)$.
}
\bigskip 

\noindent
%%%%%%%%%%%%%%%%%%%%%%%%%%%%%%%%%%%%%%%%%%%%%%%%%%%%%%%%%%%
\section{Time-decay estimates with $2 \le q \le \infty$ }
In this section, we show the time-decay estimates with $2 \le q \le \infty$ 
(not assuming $L^1$-integrability to the initial perturbation), 
that is, Theorem 3.1. 
To do that, we shall obtain the time-weighted $L^q$-energy estimates 
to $\phi$ with $2 \le q< \infty$ (cf. \cite{yoshida}). 

\medskip

\noindent
{\bf Proposition 4.1.}\quad {\it
Suppose the same assumptions in Theorem 3.1. 
For any $q \in [ \, 2, \infty)$, 
there exist positive constants $\alpha$ and $C_{\alpha,p,q}$, 
such that the unique global solution in time $\phi$ of the Cauchy problem {\rm(3.3)} 
satisfies the following $L^q$-energy estimate 
 \begin{align}
 \begin{aligned}
 &(1+t)^\alpha  \| \,\phi(t) \,\|_{L^q}^q 
  + \int ^t_0 (1+\tau )^\alpha  \int _{-\infty}^{\infty} 
    | \, \phi \, |^{q} \, \partial _x U^r 
    \, \mathrm{d}x \mathrm{d}\tau \\
 & \quad 
  + \int ^t_0 (1+\tau )^\alpha  \int _{-\infty}^{\infty} 
    | \, \phi \, |^{q-2} 
    \bigl( \, \partial _x \phi \, \bigr)^2 
    \left( \, 
    \bigl| \partial _x \phi \bigr|^{p-1} 
    + \bigl| \partial _x U^r \bigr|^{p-1}  
    \, \right) \, \mathrm{d}x \mathrm{d}\tau  \\
 & \quad 
  + \int ^t_0 (1+\tau )^\alpha  \int _{-\infty}^{\infty} 
    | \, \phi \, |^{q-2} \, 
    \left| \, 
    \bigl| \partial _x \phi + \partial _x U^r \bigr|^{p-1} 
    - \bigl| \partial _x U^r \bigr|^{p-1}  
    \, \right| \\
 & \qquad \qquad \qquad \qquad \quad \quad \quad \; \, 
   \times 
    \left| \, 
    \bigl( \, \partial _x \phi + \partial _x U^r \, \bigr)^2 
    - \bigl( \, \partial _x U^r \, \bigr)^2 \, \right| 
    \, \mathrm{d}x \mathrm{d}\tau  \\
 &\leq C_{\alpha,p,q} \| \,\phi_0 \,\|_{L^q}^q 
       + C\left( \, \alpha, \, p, \, q, \, \phi_0 \, \right) \, 
       (1+t)^{\alpha - \frac{q-2}{3p+1}}
       \quad \bigl( t \ge 0 \bigr). 
 \end{aligned}
 \end{align}
}

\medskip

The proof of Proposition 4.1 is provided by the
following two lemmas.

\medskip

\noindent
{\bf Lemma 4.1.}\quad {\it
For any $2\leq q < \infty$, 
there exist positive constants $\alpha$ and $C_q$ such that 
\begin{align}
\begin{aligned}
&(1+t)^\alpha  \| \,\phi(t) \,\|_{L^q}^q 
  + q \, (q-1)
    \int ^t_0 (1+\tau )^\alpha  \\
& \qquad \quad \quad \; \; \; \; \: \: \: \, \, 
    \times 
    \int _{-\infty}^{\infty} \int _{0}^{\phi} 
    \left( \lambda(\tilde{U}+\eta)-\lambda(\tilde{U}) \right)
    | \, \eta \, |^{q-2} \, \mathrm{d}\eta \, \bigl( \, \partial _x U^r \, \bigr) 
    \, \mathrm{d}x \mathrm{d}\tau \\
 & 
  + C_{q} 
    \int ^t_0 (1+\tau )^\alpha 
    \int _{-\infty}^{\infty} 
    | \, \phi \, |^{q-2} 
    \bigl( \, \partial _x \phi \, \bigr)^2 
    \left( \, 
    \bigl| \partial _x \phi \bigr|^{p-1} 
    + \bigl| \partial _x U^r \bigr|^{p-1}  
    \, \right) \, \mathrm{d}x \mathrm{d}\tau  \\
 & 
  + C_{q} 
    \int ^t_0 (1+\tau )^\alpha  \int _{-\infty}^{\infty} 
    | \, \phi \, |^{q-2} \, 
    \left| \, 
    \bigl| \partial _x \phi + \partial _x U^r \bigr|^{p-1} 
    - \bigl| \partial _x U^r \bigr|^{p-1}  
    \, \right| \\
 & \qquad \qquad \qquad \qquad \quad \quad \quad \; \, 
   \times 
    \left| \, 
    \bigl( \, \partial _x \phi + \partial _x U^r \, \bigr)^2 
    - \bigl( \, \partial _x U^r \, \bigr)^2 \, \right| 
    \, \mathrm{d}x \mathrm{d}\tau  \\
&\leq \| \,\phi_0 \,\|_{L^q}^q 
      + \alpha 
      \int ^t_0 (1+\tau )^{\alpha -1} 
      \| \,\phi(\tau ) \,\|_{L^p}^p 
      \, \mathrm{d}\tau \\
&\qquad \qquad \; \; \,  
 + \mu \int ^t_0 (1+\tau )^\alpha 
   \| \,\phi(\tau) \,\|_{L^\infty}^{p-1} \\
& \qquad \qquad \quad \quad \quad \quad \quad \; \, 
   \times 
   \left| \left| \, 
    \partial _x 
    \bigl( \, \bigl| \partial _x U^r \bigr|^{p-1} 
    \partial _x U^r \, \bigr) (\tau ) \, \right| \right|_{L^1} 
    \, \mathrm{d}\tau
   \quad (t \ge 0). 
\end{aligned}
\end{align}
}

\bigskip

\noindent
{\bf Lemma 4.2.}\quad {\it
Assume $p>1$ and $2\leq q < \infty$. 
We have the following interpolation inequalities.  

 \noindent
 {\rm (1)}\ \ For any $2\leq r < \infty$, there exists a positive 
 constant $C_{p,q,r}$  such that 
 \begin{align*}
 \begin{aligned}
 &\Vert \, \phi (t) \, \Vert _{L^r }
  \leq 
  C_{p,q,r} \left( \, \int _{-\infty}^{\infty} | \, \phi \, |^2 \, \mathrm{d}x \, \right)
  ^{\frac{pr+p+q-1}{(3p+q-1)r}} \\
 & \qquad \quad \quad \quad \quad \; \: \, 
         \times 
         \left( \, \int _{-\infty}^{\infty} | \, \phi \, |^{q-2} 
         \bigl| \, \partial _x \phi \, \bigr|^{p+1} \, \mathrm{d}x \, \right)
         ^{\frac{r-2}{(3p+q-1)r}} \quad \bigl( t \ge 0 \bigr). 
 \end{aligned}
 \end{align*}

 \noindent
 {\rm (2)}\ \ 
 There exists a positive 
 constant $C_{p,q}$ 
 such that 
 \begin{align*}
 \begin{aligned}
 &\Vert \, \phi (t)  \, \Vert _{L^\infty }
  \leq 
  C_{p,q} \left( \, \int _{-\infty}^{\infty} | \, \phi \, |^2 \, \mathrm{d}x \, \right)
  ^{\frac{p}{3p+q-1}} \\
 & \qquad \quad \quad \quad \quad \; \, 
         \times 
         \left( \, \int _{-\infty}^{\infty} | \, \phi \, |^{q-2} 
         \bigl| \, \partial _x \phi \, \bigr|^{p+1} \, \mathrm{d}x \, \right)
         ^{\frac{1}{3p+q-1}} \quad \bigl( t \ge 0 \bigr). 
 \end{aligned}
 \end{align*}
}

\medskip

In what follows, we first prove Lemma 4.1 and Lemma 4.2, 
and finally give the proof of Proposition 4.1. 

\medskip

{\bf Proof of Lemma 4.1}.
Multiplying the equation in (3.3) by 
$\left|\phi \right|^{q-2} \phi$ with $2\leq q < \infty$, 
we obtain the divergence form 
 \begin{align}
 \begin{aligned}
 &\partial_t\left(\frac{1}{q} \left|\, \phi \, \right|^q \right) 
  +\partial _x 
   \biggl( \, 
   \left|\, \phi \, \right|^{q-2} \phi \, 
   \bigl( f(U^r+\phi)-f(U^r) \bigr) 
   \biggr) \\ 
 &+\partial _x \left( 
   -(q-1)\int _0^{\phi} 
   \bigl( f(U^r+\eta)-f(U^r) \bigr)\left|
   \, \eta \, \right|^{q-2}
   \, \mathrm{d}\eta 
   \, \right) \\
 &+\partial _x \biggl( \, 
   -\mu \, \left|\, \phi \, \right|^{q-2} \phi \, 
   \biggr.\\
 & \quad \times \biggl. 
    \Bigl( \, 
    \bigl|\, \partial _x U^r + \partial _x \phi \, \bigr|^{p-1} 
    \bigl( \, \partial _x U^r + \partial _x \phi \, \bigr) 
    - \bigl|\, \partial _x U^r \, \bigr|^{p-1} 
      \bigl( \, \partial _x U^r \, \bigr) \, 
   \Bigr)
   \, \biggr) \\
 &+(q-1)\int _0^{\phi} 
   \bigl( \lambda(U^r+\eta)-\lambda(U^r) \bigr)\left| 
   \, \eta \, \right|^{q-2}
   \, \mathrm{d}\eta \, \bigl( \, \partial _x U^r \, \bigr) \\ 
 &+\mu \, (q-1) \, \left|\, \phi \, \right|^{q-2} \partial _x \phi \, \\
 & \quad \times 
   \Bigl( \, 
   \bigl|\, \partial _x U^r + \partial _x \phi \, \bigr|^{p-1} 
   \bigl( \, \partial _x U^r + \partial _x \phi \, \bigr) 
   - \bigl|\, \partial _x U^r \, \bigr|^{p-1} 
     \bigl( \, \partial _x U^r \, \bigr) \, 
   \Bigr)\\
 &=\mu \, \left|\, \phi \, \right|^{q-2} \phi \, 
   \partial _x 
    \bigl( \, \bigl| \partial _x U^r \bigr|^{p-1} 
    \partial _x U^r \, \bigr). % ・ｽ・ｽ ・ｽ・ｽ
 \end{aligned}
 \end{align}
Integrating (4.3) with respect to $x$, we have
 \begin{align}
 \begin{aligned}
 &\frac{1}{q} \frac{\mathrm{d}}{\mathrm{d}t} 
 \,\Vert \, \phi(t) \, \Vert_{L^q}^q \\
 &+\int ^{\infty }_{-\infty } (q-1) \int _0^{\phi} 
 \left( \lambda(U^r+\eta)-\lambda(U^r) \right)\left| 
 \, \eta \, \right|^{q-2}
 \, \mathrm{d}\eta \, \bigl( \, \partial _x U^r \, \bigr) 
 \, \mathrm{d}x \\
 &+\mu \, (q-1) \, \int ^{\infty }_{-\infty }
   \left|\, \phi \, \right|^{q-2} \partial _x \phi \, \\
 & \quad \times 
   \Bigl( \, 
   \bigl|\, \partial _x U^r + \partial _x \phi \, \bigr|^{p-1} 
   \bigl( \, \partial _x U^r + \partial _x \phi \, \bigr) 
   - \bigl|\, \partial _x U^r \, \bigr|^{p-1} 
     \bigl( \, \partial _x U^r \, \bigr) \, 
   \Bigr) 
   \, \mathrm{d}x \\
 &= \mu 
    \int ^{\infty}_{-\infty} \left|\, \phi \, \right|^{q-2} \phi 
    \, \partial _x 
    \bigl( \, \bigl| \partial _x U^r \bigr|^{p-1} 
    \partial _x U^r \, \bigr) \, \mathrm{d}x. 
 \end{aligned}
 \end{align}
By using the uniform boundedness, Lemma 3.1, 
and the following absolute equality 
with $p>1$, 
for any $a,b \in \mathbb{R}$, 
\begin{align}
\begin{aligned}
&\left(\, |\, a\, |^{p-1}a - |\, b\, |^{p-1}b \, \right) 
 \left(\, a-b \, \right) \\
&= \frac{1}{2} \, 
   \left(\, |\, a\, |^{p-1} + |\, b\, |^{p-1} \, \right) 
   \left(\, a-b \, \right)^2 
   + \frac{1}{2} \, 
     \left(\, |\, a\, |^{p-1} - |\, b\, |^{p-1} \, \right) 
     \left(\, a^2-b^2 \, \right), 
\end{aligned}
\end{align}
we have 
 \begin{align}
 \begin{aligned}
 &\frac{1}{q} \frac{\mathrm{d}}{\mathrm{d}t} 
 \,\Vert \, \phi(t) \, \Vert_{L^q}^q 
 + (q-1) \, 
   \left( \, 
    \displaystyle{\min _{|u| \leq \widetilde{C}} {\lambda}'(u)} \, 
    \right) \, 
    \int ^{\infty }_{-\infty } 
    \left|\, \phi \, \right|^{q}
    \, \partial _x U^r 
    \, \mathrm{d}x \\
 &+\frac{\mu \, (q-1)}{2} \, 
   \int ^{\infty }_{-\infty }
   | \, \phi \, |^{q-2} 
    \bigl( \, \partial _x \phi \, \bigr)^2 
    \left( \, 
    \bigl| \partial _x \phi + \partial _x U^r \bigr|^{p-1} % ・ｽ・ｽ ・ｽ・ｽ
    + \bigl| \partial _x U^r \bigr|^{p-1}  
    \, \right) 
   \, \mathrm{d}x \\
 &+\frac{\mu \, (q-1)}{2} \, 
   \int ^{\infty }_{-\infty }
   | \, \phi \, |^{q-2} \, 
    \left| \, 
    \bigl| \partial _x \phi + \partial _x U^r \bigr|^{p-1} 
    - \bigl| \partial _x U^r \bigr|^{p-1}  
    \, \right| \\
 & \qquad \qquad \qquad \qquad \quad \; \; \, \, 
   \times 
    \left| \, 
    \bigl( \, \partial _x \phi + \partial _x U^r \, \bigr)^2 
    - \bigl( \, \partial _x U^r \, \bigr)^2 \, \right| 
    \, \mathrm{d}x \\
 &\leq \mu \, \left| \,  
       \int ^{\infty}_{-\infty} \left|\, \phi \, \right|^{q-2} \phi 
       \, \partial _x 
       \bigl( \, \bigl| \partial _x U^r \bigr|^{p-1} 
       \partial _x U^r \, \bigr) \, \mathrm{d}x 
       \, \right|. 
 \end{aligned}
 \end{align}
Thus, multiplying the inequality by 
$(1+t)^{\alpha}$ with $\alpha>0$ 
and integrating over $(0,t)$ with respect to the time, 
we complete the proof of Lemma 4.1.

\medskip

{\bf Proof of Lemma 4.2}.\ 
Noting that $\phi (t, \cdot \: ) \in L^2$ 
and $\partial_x \phi (t, \cdot \: ) \in L^{p+1}$ 
imply 
$\displaystyle{\lim _{x\rightarrow \pm \infty}\phi(t,x)}=0$ 
for $t \ge 0$, 
we have 
\begin{align}
\begin{aligned}
 | \, \phi \, |^{s} 
 &\leq s \int _{-\infty}^{\infty} 
      | \, \phi \, |^{s-1}  \, 
      \bigl| \, \partial _x \phi \, \bigr|
      \, \mathrm{d}x.\\ % ・ｽ・ｽ ・ｽ・ｽ・ｽ・ｽ
\end{aligned}
\end{align}
By the Cauchy-Schwarz inequality, we have 
\begin{align}
\begin{aligned}
 | \, \phi \, |^{s} 
 &\leq s \left( \, 
       \int _{-\infty}^{\infty} 
       \left( \, | \, \phi \, |^{s-1-\frac{q-2}{p+1}} \, \right)^{\frac{p+1}{p}}
       \, \mathrm{d}x
       \right)^{\frac{p}{p+1}} \\
 & \qquad \times 
       \left( \, 
       \int _{-\infty}^{\infty} 
       | \, \phi \, |^{q-2} 
       \bigl| \, \partial _x \phi \, \bigr|^{p+1}
       \, \mathrm{d}x
       \right)^{\frac{1}{p+1}}.  % ・ｽ・ｽ ・ｽ・ｽ・ｽ・ｽ
\end{aligned}
\end{align}
Taking $s=\frac{3p+q-1}{p+1}$, we get 
\begin{align}
\begin{aligned}
 \Vert \, \phi \, \Vert _{L^\infty }^{\frac{3p+q-1}{p+1}}
 &\le \frac{3\, p+q-1}{p+1} 
     \left( \, 
     \int _{-\infty}^{\infty} 
     | \, \phi \, |^{2} 
     \, \mathrm{d}x
     \, \right)^{\frac{p}{p+1}} \\
 & \qquad \times 
     \left( \, 
     \int _{-\infty}^{\infty} 
     | \, \phi \, |^{q-2} 
     \bigl| \, \partial _x \phi \, \bigr|^{p+1}
     \, \mathrm{d}x
     \right)^{\frac{1}{p+1}},
\end{aligned}
\end{align}
and 
\begin{align}
\begin{aligned}
 \Vert \, \phi \, \Vert _{L^r }^r
 &\le \Vert \, \phi \, \Vert _{L^\infty }^{r-2}
      \Vert \, \phi \, \Vert _{L^2 }^2\\
 &\le \left( \frac{3\, p+q-1}{p+1} \right)^{\frac{(p+1)(r-2)}{3p+q-1}}
      \left( \, 
     \int _{-\infty}^{\infty} 
     | \, \phi \, |^{2} 
     \, \mathrm{d}x
     \, \right)^{\frac{pr+p+q-1}{3p+q-1}} \\
 & \qquad \times 
     \left( \, 
     \int _{-\infty}^{\infty} 
     | \, \phi \, |^{q-2} 
     \bigl| \, \partial _x \phi \, \bigr|^{p+1}
     \, \mathrm{d}x
     \right)^{\frac{r-2}{3p+q-1}}.
\end{aligned}
\end{align}
Thus, we complete the proof of Lemma 4.2. 

\medskip

{\bf Proof of Proposition 4.1}.\
By using Lemma 4.1 and Lemma 4.2, 
we shall estimate the second term, the third term and the fourth term 
on the right-hand side of (4.2) as follows: 
for any $\epsilon>0$, 
\begin{align}
\begin{aligned}
 &\alpha \int ^t_0 (1+\tau )^{\alpha -1} 
  \| \,\phi(\tau) \,\|_{L^q}^q 
  \, \mathrm{d}\tau \\
 &\le C_{\alpha ,p,q} \int ^t_0 (1+\tau )^{\alpha -1} 
      \left( \, 
      \int _{-\infty}^{\infty} 
      | \, \phi \, |^{q-2} 
      \bigl| \, \partial _x \phi \, \bigr|^{p+1}
      \, \mathrm{d}x
      \right)^{\frac{q-2}{3p+q-1}} \\
 &\qquad \qquad \qquad \qquad \qquad \quad \times 
      \left( \, 
      \int _{-\infty}^{\infty} 
      | \, \phi \, |^{2} 
      \, \mathrm{d}x
      \, \right)^{\frac{pq+p+q-1}{3p+q-1}}
      \, \mathrm{d}\tau \\
 &\le \int ^t_0 
      \left( \, 
      (1+\tau )^{\alpha} 
      \int _{-\infty}^{\infty} 
      | \, \phi \, |^{q-2} 
      \bigl| \, \partial _x \phi \, \bigr|^{p+1}
      \, \mathrm{d}x
      \, \right)^{\frac{q-2}{3p+q-1}} \\
 &\qquad \quad \quad \quad \quad \times 
      C_{\alpha ,p,q}\,  (1+\tau )^{\alpha - 1 -\frac{\alpha (q-2)}{3p+q-1}}
         \| \,\phi(\tau) \,\|_{L^2}^{\frac{2(pq+p+q-1)}{3p+q-1}}
         \, \mathrm{d}\tau \\
 &\le \epsilon \int ^t_0 (1+\tau )^{\alpha} 
      \left( \, 
      \int _{-\infty}^{\infty} 
      | \, \phi \, |^{q-2} 
      \bigl| \, \partial _x \phi \, \bigr|^{p+1}
      \, \mathrm{d}x
      \right) \, \mathrm{d}\tau \\
 &\quad \quad \quad \quad \quad 
         +C_{\alpha ,p,q}(\epsilon ) 
         \int ^t_0 (1+\tau )^{\alpha - {\frac{3p+q-1}{3p+1}}}
         \| \,\phi(\tau ) \,\|_{L^2}^{\frac{2(pq+p+q-1)}{3p+1}} 
         \, \mathrm{d}\tau, 
\end{aligned}
\end{align}
\begin{align}
\begin{aligned}
 &q \int ^t_0 (1+\tau )^{\alpha} 
  \| \,\phi(\tau ) \,\|_{L^\infty}^{q-1} 
  \left| \,  
       \int ^{\infty}_{-\infty} \left|\, \phi \, \right|^{q-2} \phi 
       \, \partial _x 
       \bigl( \, \bigl| \partial _x U^r \bigr|^{p-1} 
       \partial _x U^r \, \bigr) \, \mathrm{d}x 
       \, \right|
  \, \mathrm{d}\tau \\
 &\le C_{p,q} \int ^t_0 (1+\tau )^{\alpha} 
      \left( \, 
      \int _{-\infty}^{\infty} 
      | \, \phi \, |^{q-1} 
      \bigl| \, \partial _x \phi \, \bigr|^{p+1}
      \, \mathrm{d}x
      \right)^{\frac{q-1}{3p+q-1}} \\
 &\qquad \quad \quad \times 
      \left( \, 
      \int _{-\infty}^{\infty} 
      | \, \phi \, |^{2} 
      \, \mathrm{d}x
      \, \right)^{\frac{p(q-1)}{3p+q-1}}
      \left|\left| \, 
      \partial _x 
      \bigl( \, \bigl| \partial _x U^r \bigr|^{p-1} 
      \partial _x U^r \, \bigr)(\tau ) 
      \, \right| \right|_{L^1}
      \, \mathrm{d}\tau \\
 &\le \int ^t_0 
      \left( \, 
      (1+\tau )^{\alpha} 
      \int _{-\infty}^{\infty} 
      | \, \phi \, |^{q-2} 
      \bigl| \, \partial _x \phi \, \bigr|^{p+1}
      \, \mathrm{d}x
      \, \right)^{\frac{q-2}{3p+q-1}} \\
 &\quad \: \: \: \times 
         C_{p,q}\,  (1+\tau )^{\alpha -\frac{\alpha (q-1)}{3p+q-1}} 
% &\quad \: \: \: \times 
         \| \,\phi(\tau) \,\|_{L^2}^{\frac{2p(q-1)}{3p+q-1}} 
         \left|\left| \, 
         \partial _x 
         \bigl( \, \bigl| \partial _x U^r \bigr|^{p-1} 
         \partial _x U^r 
         \, \bigr)(\tau ) \, \right| \right|_{L^1}
         \, \mathrm{d}\tau \\
 &\le \epsilon \int ^t_0 (1+\tau )^{\alpha} 
      \left( \, 
      \int _{-\infty}^{\infty} 
      | \, \phi \, |^{q-1} 
      \bigl| \, \partial _x \phi \, \bigr|^{p+1}
      \, \mathrm{d}x
      \right) 
      \, \mathrm{d}\tau \\
 &\quad 
         + C_{p,q}(\epsilon) \int ^t_0 
          (1+\tau )^{\alpha} \| \,\phi(\tau ) \,\|
          _{L^2}^{{\frac{2}{3}}(q-1)}
          \left|\left| \, 
          \partial _x 
          \bigl( \, \bigl| \partial _x U^r \bigr|^{p-1} 
          \partial _x U^r \, \bigr)(\tau ) \, \right| \right|_{L^1}
          ^{\frac{3p+q-1}{3p}} 
          \, \mathrm{d}\tau. 
\end{aligned}
\end{align}
Substituting (4.11) and (4.12) into (4.2), 
we have 
\begin{align}
\begin{aligned}
 &(1+t)^\alpha  \| \,\phi(t) \,\|_{L^q}^q 
  + \int ^t_0 (1+\tau )^\alpha  \int _{-\infty}^{\infty} 
    | \, \phi \, |^{q} \, \partial _x U^r 
    \, \mathrm{d}x \mathrm{d}\tau \\
 & %\quad 
  + \int ^t_0 (1+\tau )^\alpha  \int _{-\infty}^{\infty} 
    | \, \phi \, |^{q-2} 
    \bigl( \, \partial _x \phi \, \bigr)^2 
    \left( \, 
    \bigl| \partial _x \phi \bigr|^{p-1} 
    + \bigl| \partial _x U^r \bigr|^{p-1}  
    \, \right) \, \mathrm{d}x \mathrm{d}\tau  \\
 & %\quad 
  + \int ^t_0 (1+\tau )^\alpha  \int _{-\infty}^{\infty} 
    | \, \phi \, |^{q-2} \, 
    \left| \, 
    \bigl| \partial _x \phi + \partial _x U^r \bigr|^{p-1} 
    - \bigl| \partial _x U^r \bigr|^{p-1}  
    \, \right| \\
 & \qquad \qquad \qquad \qquad \quad \quad \; \, \, 
   \times 
    \left| \, 
    \bigl( \, \partial _x \phi + \partial _x U^r \, \bigr)^2 
    - \bigl( \, \partial _x U^r \, \bigr)^2 \, \right| 
    \, \mathrm{d}x \mathrm{d}\tau  \\
 &\leq C_{\alpha,p,q} \| \,\phi_0 \,\|_{L^q}^q 
       + C_{\alpha,p,q} 
         \int ^t_0 (1+\tau )^{\alpha - {\frac{3p+q-1}{3p+1}}}
         \| \,\phi(\tau ) \,\|_{L^2}^{\frac{2(pq+p+q-1)}{3p+1}} 
         \, \mathrm{d}\tau \\
 & \qquad \qquad \qquad \quad \: 
   + C_{p,q} 
     \int ^t_0 
     (1+\tau )^{\alpha} \| \,\phi(\tau ) \,\|
     _{L^2}^{{\frac{2}{3}}(q-1)} \\
 & \qquad \qquad \qquad \qquad \quad \quad \times 
     \left|\left| \, 
     \partial _x 
     \bigl( \, \bigl| \partial _x U^r \bigr|^{p-1} 
     \partial _x U^r \, \bigr)(\tau ) \, \right| \right|_{L^1}
     ^{\frac{3p+q-1}{3p}} 
     \, \mathrm{d}\tau.
%     \quad \bigl( t \ge T_0 \bigr). 
\end{aligned}
\end{align}
By using the $L^2$-boundedness of $\phi$, (3.6), 
and 
\begin{align}
\begin{aligned}
\left|\left| \, 
\partial _x 
\bigl( \, \bigl| \partial _x U^r \bigr|^{p-1} 
\partial _x U^r \, \bigr)(\tau ) \, \right| \right|_{L^1}
& \leq p \, \| \, \partial _x U^r (t) \,\|_{L^{\infty}}^{p-1} 
       \| \, \partial _x^2 U^r (t) \,\|_{L^1} \\
& \leq C_{p} (1+t)^{-p}, 
\end{aligned}
\end{align}
we estimate the each terms on the right-hand side of (4.13) 
as follows: 
\begin{align}
\begin{aligned}
&C_{\alpha ,p,q} \, 
         \int ^t_0 (1+\tau )^{\alpha - {\frac{3p+q-1}{3p+1}}} 
         \| \,\phi(\tau ) \,\|_{L^2}^{\frac{2(pq+p+q-1)}{3p+1}} 
         \, \mathrm{d}\tau \\
&\leq C_{\alpha ,p,q} 
      \bigl( \, C_{p}(\phi_0) \, \bigr)^{\frac{pq+p+q-1}{3p+1}} 
      \int ^t_0 (1+\tau )^{\alpha - {\frac{3p+q-1}{3p+1}}} 
      \, \mathrm{d}\tau \\
&\leq C_{\alpha ,p,q} 
      \bigl( \, C_{p}(\phi_0) \, \bigr)^{\frac{pq+p+q-1}{3p+1}} 
      (1+t )^{\alpha - {\frac{q-2}{3p+1}}}, 
\end{aligned}
\end{align}
\begin{align}
\begin{aligned}
&C_{p,q} \int ^t_0 
          (1+\tau )^{\alpha} \| \,\phi(\tau ) \,\|
          _{L^2}^{{\frac{2}{3}}(q-1)}
          \left|\left| \, 
          \partial _x 
          \bigl( \, \bigl| \partial _x U^r \bigr|^{p-1} 
          \partial _x U^r \, \bigr)(\tau ) \, \right| \right|_{L^1}
          ^{\frac{3p+q-1}{3p}} 
          \, \mathrm{d}\tau \\
&\leq C_{p,q} 
      \bigl( \, C_{p}(\phi_0) \, \bigr)^{{\frac{1}{3}}(q-1)} 
      \int ^t_0 (1+\tau )^{\alpha - {\frac{3p+q-1}{3}} } 
      \, \mathrm{d}\tau \\
&\leq C_{p,q} 
      \bigl( \, C_{p}(\phi_0) \, \bigr)^{{\frac{1}{3}}(q-1)} 
      (1+t )^{\alpha - {\frac{3p+q-4}{3}}}. 
\end{aligned}
\end{align}
Substituting (4.15) and (4.16) into (4.13), 
we get (4.1). 
Thus the proof of Proposition 4.1 is complete. 
In particular, it follows that 
\begin{align}
\begin{aligned}
\| \,\phi(t) \,\|_{L^q} 
\leq C( \, p, \,q, \, \phi_0 \, ) \, 
     (1+t)^{-{\frac{1}{3p+1}} \left( 1 - {\frac{2}{q}} \right)}
\end{aligned}
\end{align}
for $2 \leq q < \infty$. 

\medskip

{\bf Proof of Theorem 3.1.}\ 
We already have proved the decay estimate 
of $\| \,\phi(t) \,\|_{L^q}$ with $2 \le q < \infty$. 
Therefore we only 
show the $L^{\infty}$-estimate. 
We first note by Lemma 2.2 that 
\begin{align}                
\begin{aligned}
&\bigl|\bigl| 
 \, \partial _x \phi (t)\, 
 \bigr|\bigr| _{L^{p+1}}^{p+1}\\
&\leq \bigl|\bigl| 
      \, \partial _x u (t)\, 
      \bigr|\bigr| _{L^{p+1}}^{p+1} 
      + \bigl|\bigl| 
        \, \partial _x U^r (t)\, 
        \bigr|\bigr| _{L^{p+1}}^{p+1} \\
&\leq C\bigl( \, \| \, \phi_0 \, \|_{L^2}, \: 
       \| \, \partial _x u_0 \, \|_{L^{p+1}} 
       \, \bigr) 
       + C_{p}(1+t)^{-p}. 
\end{aligned}
\end{align}
We use the following Gagliardo-Nirenberg inequality: 
\begin{align}
\| \,\phi(t) \,\|_{L^\infty }
\le C_{q,\theta}
\| \,\phi(t) \,\|_{L^q}^{1-\theta}
\| \,\partial _x \phi(t) \,\|_{L^{p+1}}^{\theta}
\end{align}
for any 
$(q,\theta)\in [\, 1,\infty)\times (0,1\, ]$ satisfying 
$$
%\, ;\, 
\frac{p}{p+1} \, {\theta}=(1-\theta) \, \frac{1}{q}. 
$$
Substituting (4.17) and (4.18) into (4.19), 
we have 
\begin{align}
\begin{aligned}
\| \,\phi(t) \,\|_{L^\infty }
&\le C( \, p, \,q, \, \theta, \, \phi_0, \, \partial _x u_0 \, ) \, 
     (1+t)^{-{\frac{1}{3p+1}} \left( 1 - {\frac{2}{q}} \right)(1-\theta)}\\
&\le C( \, p, \, \theta, \, \phi_0, \, \partial _x u_0 \, ) \, 
     (1+t)^{-{\frac{1}{3p+1}} + {\frac{\theta}{p+1}}}
\end{aligned}
\end{align}
for $\theta \in (0,1\, ]$. 
Consequently, we do complete the proof of Theorem 3.1.
\bigskip 

\noindent
%%%%%%%%%%%%%%%%%%%%%%%%%%%%%%%%%%%%%%%%%%%%%%%%%%%%%%%%%%%%%%%%%%%%%%%%%%%%%
\section{Time-decay estimates with $1 \le q \le \infty$ }
In this section, 
we show the time-decay estimates 
with $1 \le q \le \infty$ and 
time-decay estimate for the higher order derivative 
in the $L^{p+1}$-norm, 
in the case where $\phi_0 \in L^1 \cap L^2$ 
with $\partial _x u_0 \in L^{p+1}$, 
that is, Theorem 3.2. 
Then, we 
first establish the $L^1$-estimate 
to the solution $\phi$ of the reformulated Cauchy problem (3.3). 
To do that, we use the Friedrichs mollifier $\rho_\delta \ast $, 
where, 
$\rho_\delta(\phi):=\frac{1}{\delta}\rho \left( \frac{\phi}{\delta}\right)$
with 
\begin{align*}
\begin{aligned}
&\rho \in C^{\infty}_0(\mathbb{R}), \quad
\rho (\phi)\geq 0 \quad (\phi \in \mathbb{R}),
\\
&\mathrm{supp} \{\rho \} \subset 
\left\{\phi \in \mathbb{R}\: \left|\:  |\,\phi \, |\le 1 \right. \right\},\quad  
\int ^{\infty}_{-\infty} \rho (\phi)\, \mathrm{d}\phi=1.
\end{aligned}
\end{align*}

%We record some useful properties as 
Some useful properties of the mollifier are as follows. 
\medskip

\noindent
{\bf Lemma 5.1.}\quad{\it
%For each $\phi \in \mathbb{R}$, we have the following:

\noindent
{\rm (i)}\ 
$\displaystyle{\lim_{\delta \to 0}\, 
\left( \rho_\delta \ast \mathrm{sgn} \right)(\phi)}
= \mathrm{sgn} (\phi)
\qquad (\phi \in \mathbb{R}),$

\noindent
{\rm (ii)}\ 
$\displaystyle{\lim_{\delta \to 0}\, 
\int^{\phi}_0\left( \rho_\delta \ast \mathrm{sgn} \right)(\eta)
\, \mathrm{d}\eta}
=|\, \phi \, |
\qquad (\phi \in \mathbb{R}),$

\medskip

\noindent
{\rm (iii)}\ 
$\Big. \left( \rho_\delta \ast \mathrm{sgn} \right) \Bigr|_{\phi=0} =0,$

\medskip

\noindent
{\rm (iv)}\ 
$\displaystyle{ \frac{\mathrm{d}}{\mathrm{d}\phi}\, 
\left( \rho_\delta \ast \mathrm{sgn} \right)(\phi)}
=2\, \rho_\delta(\phi)
\ge 0
\qquad (\phi \in \mathbb{R}),$

\medskip

\noindent
where 
$$
\left( \rho_\delta \ast \mathrm{sgn} \right)(\phi)
:=\int^{\infty}_{-\infty}
\rho_\delta(\phi-y)\, \mathrm{sgn}(y)\, \mathrm{d}y
\qquad (\phi \in \mathbb{R})
$$
and 
\begin{equation*}
\mathrm{sgn}(\phi):=\left\{
\begin{array}{ll}
-1& \quad \; \bigl(\, \phi < 0 \, \bigr),\\[3pt]
\, \, 0& \quad \; \bigl(\, \phi = 0 \, \bigr),\\[3pt]
\, \, 1& \quad \; \bigl(\, \phi > 0 \, \bigr).
\end{array}
\right.
\end{equation*} 
}

\noindent
Making use of Lemma 5.1, 
we obtain the following $L^1$-estimate. 

\medskip

\noindent
{\bf Proposition 5.1.}\quad {\it
Assume that the same assumptions in Theorem 3.2. 
For any $p>1$, 
the unique global solution in time $\phi$ 
of the Cauchy problem {\rm(3.3)} 
satisfies the following $L^1$-estimate 
\begin{align}
\begin{aligned}
\| \, \phi(t) \, \|_{L^1}
\le \| \,\phi_0 \,\|_{L^1} + C_{p}(1+t)^{-(p-1)}
\quad \bigl( t \ge 0 \bigr). 
\end{aligned}
\end{align}
}

\medskip

\noindent
{\bf Proof of Proposition 5.1.}\quad
Multiplying the equation in the problem (3.3) by 
$\left( \rho_\delta \ast \mathrm{sgn} \right)(\phi)$, 
we obtain the divergence form 
\medskip
\begin{align}
\begin{aligned}
&\partial_t\left( 
\int^{\phi}_0 
\left( \rho_\delta \ast \mathrm{sgn} \right)(\eta)
\, \mathrm{d}\eta 
\right) \\
&+\partial _x 
  \biggl( \, 
  \left( \rho_\delta \ast \mathrm{sgn} \right)(\phi) \, 
  \bigl( f(U^r+\phi)-f(U^r) \bigr) 
  \biggr) \\ 
&+\partial _x \left( 
  -\int _0^{\phi} 
  \bigl( f(U^r+\eta)-f(U^r) \bigr)
  \frac{\mathrm{d}\left( \rho_\delta \ast \mathrm{sgn} \right)}
  {\mathrm{d}\phi}
  (\eta)
  \, \mathrm{d}\eta 
  \, \right) \\
&+\partial _x \biggl( \, 
  -\mu \, \left( \rho_\delta \ast \mathrm{sgn} \right)(\phi) \biggr. \\
& \quad \times \biggl. 
   \Bigl( \, 
    \bigl|\, \partial _x U^r + \partial _x \phi \, \bigr|^{p-1} 
    \bigl( \, \partial _x U^r + \partial _x \phi \, \bigr) 
    - \bigl|\, \partial _x U^r \, \bigr|^{p-1} 
      \bigl( \, \partial _x U^r \, \bigr) \, 
   \Bigr)
  \, \biggr) \\
&+\int _0^{\phi} 
  \bigl( \lambda(U^r+\eta)-\lambda(U^r) \bigr) 
  \frac{\mathrm{d}\left( \rho_\delta \ast \mathrm{sgn} \right)}
  {\mathrm{d}\phi}(\eta) 
  \, \mathrm{d}\eta \, \bigl( \, \partial _x U^r \, \bigr) \\ 
&+\mu \, 
  \frac{\mathrm{d}\left( \rho_\delta \ast \mathrm{sgn} \right)}
  {\mathrm{d}\phi}(\phi) \, 
  \partial _x \phi \\
& \quad \times 
  \Bigl( \, 
    \bigl|\, \partial _x U^r + \partial _x \phi \, \bigr|^{p-1} 
    \bigl( \, \partial _x U^r + \partial _x \phi \, \bigr) 
    - \bigl|\, \partial _x U^r \, \bigr|^{p-1} 
      \bigl( \, \partial _x U^r \, \bigr) \, 
   \Bigr) \\
& =\mu \, 
   \left( \rho_\delta \ast \mathrm{sgn} \right)(\phi) \, 
   \partial _x 
   \bigl( \, \bigl| \partial _x U^r \bigr|^{p-1} 
   \partial _x U^r \, \bigr). 
\end{aligned}
\end{align}
By using (4.5) and 
integrating (5.2) with respect to $x$ and $t$, we have
\begin{align}
\begin{aligned}
&\int ^{\infty}_{-\infty} \int^{\phi (t) }_0 
 \left( \rho_\delta \ast \mathrm{sgn} \right)(\eta)
 \, \mathrm{d}\eta \, \mathrm{d}x \\
&+ \int _0^t \int ^{\infty}_{-\infty} \int^{\phi}_0
   \bigl( \lambda(U^r+\eta)-\lambda(U^r) \bigr) 
   \frac{\mathrm{d}\left( \rho_\delta \ast \mathrm{sgn} \right)}
   {\mathrm{d}\phi}(\eta)
   \, \mathrm{d}\eta \, \bigl( \, \partial _x U^r \, \bigr)
   \, \mathrm{d}x \mathrm{d}\tau \\
&+ \frac{\mu}{2} \int _0^t \int ^{\infty}_{-\infty} 
   \frac{\mathrm{d}\left( \rho_\delta \ast \mathrm{sgn} \right)}
   {\mathrm{d}\phi}(\phi) \, 
   \bigl( \, \partial _x \phi \, \bigr)^2 
    \left( \, 
    \bigl| \partial _x \phi + \partial _x U^r \bigr|^{p-1} % ・ｽ・ｽ ・ｽ・ｽ
    + \bigl| \partial _x U^r \bigr|^{p-1}  
    \, \right)
   \, \mathrm{d}x \mathrm{d}\tau \\
&+ \frac{\mu}{2} \int _0^t \int ^{\infty}_{-\infty} 
   \frac{\mathrm{d}\left( \rho_\delta \ast \mathrm{sgn} \right)}
   {\mathrm{d}\phi}(\phi) \, 
   \left| \, 
    \bigl| \partial _x \phi + \partial _x U^r \bigr|^{p-1} 
    - \bigl| \partial _x U^r \bigr|^{p-1}  
    \, \right| \\
& \qquad \qquad \qquad \qquad \quad \quad \quad \; \, \, 
   \times 
    \left| \, 
    \bigl( \, \partial _x \phi + \partial _x U^r \, \bigr)^2 
    - \bigl( \, \partial _x U^r \, \bigr)^2 \, \right| 
   \, \mathrm{d}x \mathrm{d}\tau \\
&= \int ^{\infty}_{-\infty} \int^{\phi _0}_0 
   \left( \rho_\delta \ast \mathrm{sgn} \right)(\eta)
   \, \mathrm{d}\phi \, \mathrm{d}x \\
&\qquad \qquad  \qquad \; 
 + \mu 
   \int _0^t \int ^{\infty}_{-\infty} 
   \left( \rho_\delta \ast \mathrm{sgn} \right)(\phi) \, 
   \partial _x 
   \bigl( \, \bigl| \partial _x U^r \bigr|^{p-1} 
   \partial _x U^r \, \bigr)
   \, \mathrm{d}x \mathrm{d}\tau .
\end{aligned}
\end{align}
%We note 
By using Lemma 5.1, we note that for $t\in [\, 0,\infty)$,
\begin{align}
\begin{aligned}
& \int ^{\infty}_{-\infty} \int^{\phi (t)}_0 
  \bigl( \lambda(U^r+\eta)-\lambda(U^r) \bigr) 
  \frac{\mathrm{d}\left( \rho_\delta \ast \mathrm{sgn} \right)}
  {\mathrm{d}\phi}(\eta)
  \, \mathrm{d}\eta \, \bigl( \, \partial _x U^r \, \bigr)
  \, \mathrm{d}x \\
& \ge 2 \, 
      \left( \, 
      \displaystyle{\min _{|u| \leq \widetilde{C}} {\lambda}'(u)} \, 
      \right) \, 
      \int ^{\infty }_{-\infty } 
      \left|\, 
      \int^{| \phi (t) |}_0 
      \eta \, \rho_\delta (\eta ) 
      \mathrm{d}\eta 
      \, \right|
      \, \bigl( \, \partial _x U^r \, \bigr)
      \, \mathrm{d}x 
      \ge 0, 
\end{aligned}
\end{align}
\begin{align}
\left| \, 
\int^{\phi}_0\left( \rho_\delta \ast \mathrm{sgn} \right)(\eta)
\, \mathrm{d}\eta 
\, \right|
\le \left( \rho_\delta \ast \mathrm{sgn} \right)
    (|\,  \phi \,|)\, |\,  \phi \,|
\le |\,  \phi \,|, 
\end{align}
\begin{align}
\displaystyle {\lim_{\delta\rightarrow 0}}
\int ^{\infty}_{-\infty} \int^{\phi (t) }_0
\left( \rho_\delta \ast \mathrm{sgn} \right)(\eta)
\, \mathrm{d}\eta 
= \| \,\phi(t) \,\|_{L^1}, 
\end{align}
and we get 
\begin{align}
\begin{aligned}
&\| \,\phi (t) \,\|_{L^1} 
 \le \| \,\phi_0 \,\|_{L^1} \\
&\qquad \quad 
   + \displaystyle {\mu \, \lim_{\delta\rightarrow 0}}
   \int^t_0 \left|\, \int ^{\infty}_{-\infty} 
   \left( \rho_\delta \ast \mathrm{sgn} \right)(\phi) \, 
   \partial _x 
   \bigl( \, \bigl| \partial _x U^r \bigr|^{p-1} 
   \partial _x U^r \, \bigr) 
   \, \mathrm{d}x \, \right|\, \mathrm{d}\tau. 
\end{aligned}
\end{align}
By using (4.14), we can also get 
\begin{align}
\begin{aligned}
& \displaystyle {\lim_{\delta\rightarrow 0}}
  \left|\, \int ^{\infty}_{-\infty} 
  \left( \rho_\delta \ast \mathrm{sgn} \right)(\phi) \, 
  \partial _x 
  \bigl( \, \bigl| \partial _x U^r \bigr|^{p-1} 
  \partial _x U^r \, \bigr) 
  \, \mathrm{d}x \, \right|\, (t) \\
& \leq  \left(\, \int ^{\infty}_{-\infty} 
  \bigl| \, \mathrm{sgn} (\phi) \, \bigr| \, 
  \Bigl| \, 
  \partial _x 
  \bigl( \, \bigl| \partial _x U^r \bigr|^{p-1} 
  \partial _x U^r \, \bigr) 
  \, \Bigr| 
  \, \mathrm{d}x \, \right)\, (t) \\
& \leq C_{p}(1+t)^{-p}
  \quad \bigl( t \ge 0 \bigr). 
\end{aligned}
\end{align}
Then, substituting (5.8) into (5.7), 
we have the desired $L^1$-estimate (5.1). 

\medskip
Next, we show %have 
the time-weighted $L^q$-energy estimates to $\phi$. 

\medskip

\noindent
{\bf Proposition 5.2.}\quad {\it
Suppose the same assumptions in Theorem 3.2. 
For any $q \in [ \, 1, \infty)$, 
there exist positive constants $\alpha$ and $C_{\alpha,p,q}$, 
such that the unique global solution in time $\phi$ of the Cauchy problem {\rm(3.3)} 
satisfies the following $L^q$-energy estimate 
 \begin{align}
 \begin{aligned}
 &(1+t)^\alpha  \| \,\phi(t) \,\|_{L^q}^q 
  + \int ^t_0 (1+\tau )^\alpha  \int _{-\infty}^{\infty} 
    | \, \phi \, |^{q} \, \partial _x U^r 
    \, \mathrm{d}x \mathrm{d}\tau \\
 & \quad 
  + \int ^t_0 (1+\tau )^\alpha  \int _{-\infty}^{\infty} 
    | \, \phi \, |^{q-2} 
    \bigl( \, \partial _x \phi \, \bigr)^2 
    \left( \, 
    \bigl| \partial _x \phi \bigr|^{p-1} 
    + \bigl| \partial _x U^r \bigr|^{p-1}  
    \, \right) \, \mathrm{d}x \mathrm{d}\tau  \\
 & \quad 
  + \int ^t_0 (1+\tau )^\alpha  \int _{-\infty}^{\infty} 
    | \, \phi \, |^{q-2} \, 
    \left| \, 
    \bigl| \partial _x \phi + \partial _x U^r \bigr|^{p-1} 
    - \bigl| \partial _x U^r \bigr|^{p-1}  
    \, \right| \\
 & \qquad \qquad \qquad \qquad \quad \quad \quad \; \, 
   \times 
    \left| \, 
    \bigl( \, \partial _x \phi + \partial _x U^r \, \bigr)^2 
    - \bigl( \, \partial _x U^r \, \bigr)^2 \, \right| 
    \, \mathrm{d}x \mathrm{d}\tau  \\
 &\leq C_{\alpha,p,q} \| \,\phi_0 \,\|_{L^q}^q 
       + C\left( \, \alpha, \, p, \, q, \, \phi_0 \, \right) \, 
       (1+t)^{\alpha - \frac{q-1}{2p}}
       \quad \bigl( t \ge 0 \bigr). 
 \end{aligned}
 \end{align}
}

\medskip

The proof of Proposition 5.2 is given by the
following two lemmas.

\medskip

\noindent
{\bf Lemma 5.2.}\quad {\it
For any 
$1\leq q < \infty$, 
there exist positive constants $\alpha$ and $C_q$ such that 
\begin{align}
\begin{aligned}
&(1+t)^\alpha  \| \,\phi(t) \,\|_{L^q}^q 
  + q \, (q-1)
    \int ^t_0 (1+\tau )^\alpha  \\
& \qquad \quad \quad \; \; \; \; \: \: \: \, \, 
    \times 
    \int _{-\infty}^{\infty} \int _{0}^{\phi} 
    \left( \lambda(U^r+\eta)-\lambda(U^r) \right) % ・ｽ・ｽ ・ｽ・ｽ
    | \, \eta \, |^{q-2} \, \mathrm{d}\eta \, \bigl( \, \partial _x U^r \, \bigr) 
    \, \mathrm{d}x \mathrm{d}\tau \\
 & 
  + C_{q} 
    \int ^t_0 (1+\tau )^\alpha 
    \int _{-\infty}^{\infty} 
    | \, \phi \, |^{q-2} 
    \bigl( \, \partial _x \phi \, \bigr)^2 
    \left( \, 
    \bigl| \partial _x \phi \bigr|^{p-1} 
    + \bigl| \partial _x U^r \bigr|^{p-1}  
    \, \right) \, \mathrm{d}x \mathrm{d}\tau  \\
 & 
  + C_{q} 
    \int ^t_0 (1+\tau )^\alpha  \int _{-\infty}^{\infty} 
    | \, \phi \, |^{q-2} \, 
    \left| \, 
    \bigl| \partial _x \phi + \partial _x U^r \bigr|^{p-1} 
    - \bigl| \partial _x U^r \bigr|^{p-1}  
    \, \right| \\
 & \qquad \qquad \qquad \qquad \quad \quad \quad \; \, 
   \times 
    \left| \, 
    \bigl( \, \partial _x \phi + \partial _x U^r \, \bigr)^2 
    - \bigl( \, \partial _x U^r \, \bigr)^2 \, \right| 
    \, \mathrm{d}x \mathrm{d}\tau  \\
&\leq \| \,\phi_0 \,\|_{L^q}^q 
      + \alpha 
      \int ^t_0 (1+\tau )^{\alpha -1} 
      \| \,\phi(\tau ) \,\|_{L^p}^p 
      \, \mathrm{d}\tau \\
&\qquad \qquad \; \; \,  
 + \mu \int ^t_0 (1+\tau )^\alpha 
   \| \,\phi(\tau) \,\|_{L^\infty}^{p-1} \\
& \qquad \qquad \quad \quad \quad \quad \quad \; \, 
   \times 
   \left| \left| \, 
    \partial _x 
    \bigl( \, \bigl| \partial _x U^r \bigr|^{p-1} 
    \partial _x U^r \, \bigr) (\tau ) \, \right| \right|_{L^1} 
    \, \mathrm{d}\tau
   \quad (t \ge 0). 
\end{aligned}
\end{align}
}

\bigskip

\noindent
{\bf Lemma 5.3.}\quad {\it
Assume $p>1$ and $1\leq q < \infty$. 
We have the following interpolation inequalities.  

 \noindent
 {\rm (1)}\ \ For any $1\leq r < \infty$, there exists a positive 
 constant $C_{p,q,r}$  such that 
 \begin{align*}
 \begin{aligned}
 &\Vert \, \phi (t) \, \Vert _{L^r }
  \leq 
  C_{p,q,r} \left( \, \int _{-\infty}^{\infty} | \, \phi \, | \, \mathrm{d}x \, \right)
  ^{\frac{pr+p+q-1}{(2p+q-1)r}} \\
 & \qquad \quad \quad \quad \quad \; \: \, 
         \times 
         \left( \, \int _{-\infty}^{\infty} | \, \phi \, |^{q-2} % ・ｽ・ｽ ・ｽ・ｽ・ｽ・ｽ・ｽ・ｽ
         \bigl| \, \partial _x \phi \, \bigr|^{p+1} \, \mathrm{d}x \, \right)
         ^{\frac{r-1}{(2p+q-1)r}} \quad \bigl( t \ge 0 \bigr). 
 \end{aligned}
 \end{align*}

 \noindent
 {\rm (2)}\ \ 
 There exists a positive 
 constant $C_{p,q}$ 
 such that 
 \begin{align*}
 \begin{aligned}
 &\Vert \, \phi (t)  \, \Vert _{L^\infty }
  \leq 
  C_{p,q} \left( \, \int _{-\infty}^{\infty} | \, \phi \, | \, \mathrm{d}x \, \right)
  ^{\frac{p}{2p+q-1}} \\
 & \qquad \quad \quad \quad \quad \; \, 
         \times 
         \left( \, \int _{-\infty}^{\infty} | \, \phi \, |^{q-2} % ・ｽ・ｽ ・ｽ・ｽ・ｽ・ｽ・ｽ・ｽ
         \bigl| \, \partial _x \phi \, \bigr|^{p+1} \, \mathrm{d}x \, \right)
         ^{\frac{1}{2p+q-1}} \quad \bigl( t \ge 0 \bigr). 
 \end{aligned}
 \end{align*}
}

\medskip

\noindent
The proofs of Lemma 5.2, Lemma 5.3 and Proposition 5.2 are given in 
the same way as those of Lemma 4.1, Lemma 4.2 and Proposition 4.1, 
so we omit them. 
We particularly note that we have by Proposition 5.2 
\begin{equation}
\| \,\phi (t) \,\|_{L^q}
\leq C( \, p, \, q, \, \phi_0 \, )\, (1+t)^
     {-\frac{1}{2p}\left( 1- \frac{1}{q} \right)} 
\end{equation}
for $1\leq q < \infty$. 

\bigskip

We shall finally obtain the time-decay estimates 
for the higher order derivatives, that is, 
$\partial _x \phi$ and $\partial _x u$, 
and also get the $L^{\infty}$-estimate for $\phi$. 

\medskip

\noindent
{\bf Proposition 5.3.}\quad {\it
Suppose the same assumptions in Theorem 3.2. 
%For any $q \in [ \, 1, \infty)$, 
There exist positive constants $\alpha$ and $C_{\alpha,p}$, 
such that the unique global solution in time $\phi$ of the Cauchy problem {\rm(3.3)} 
satisfies the following $L^{p+1}$-energy estimate 
\begin{align}
\begin{aligned}
 &(1+t)^\alpha  
  \| \, \partial _x u (t)\, \| _{L^{p+1}}^{p+1} 
  + %\mu \, p^2 \, (p+1) 
    \int ^t_0 (1+\tau )^\alpha  \int _{-\infty}^{\infty} 
    \bigl| \, \partial _x u \, \bigr|^{2(p-1)} 
    \left( \, \partial _x^2 u \, \right)^2 
    \, \mathrm{d}x \mathrm{d}\tau \\
 & \qquad \qquad \qquad \quad \quad \quad \, 
  + \int ^t_0 (1+\tau )^\alpha  
    \| \, \partial _x u (\tau) \, \| _{L^{p+2}}^{p+2} 
    \, \mathrm{d}\tau  \\
 &\leq C_{\alpha,p} 
       \| \, \partial _x u_{0} \, \| _{L^{p+1}}^{p+1} 
       + C\left( \, \alpha, \, p, \, \phi_0, \, \partial _x u_{0} \, \right) \, 
       (1+t)^{\alpha - \frac{1}{2p}}
       \quad \bigl( t \ge 0 \bigr). 
\end{aligned}
\end{align}
}

\medskip

\noindent
To obtain Proposition 5.3, we first show the following. 

\medskip

\noindent
{\bf Lemma 5.4.}\quad {\it
It follows that 
\begin{align}
\begin{aligned}
 &(1+t)^\alpha  
  \| \, \partial _x u (t)\, \| _{L^{p+1}}^{p+1} \\
 & \quad 
  + \mu \, p^2 \, (p+1) 
    \int ^t_0 (1+\tau )^\alpha  \int _{-\infty}^{\infty} 
    \bigl| \, \partial _x u \, \bigr|^{2(p-1)} 
    \left( \, \partial _x^2 u \, \right)^2 
    \, \mathrm{d}x \mathrm{d}\tau \\
 & \qquad \quad \quad \; \; \; \: \, 
  + p \int ^t_0 (1+\tau )^\alpha  
    \int _{\partial _x u \geq 0} 
    f''(u) \left|\, \partial _x u \, \right|^{p+2} 
    \, \mathrm{d}x \mathrm{d}\tau  \\
 &= \| \, \partial _x u_{0} \, \| _{L^{p+1}}^{p+1}
    + \alpha \int ^t_0 (1+\tau )^{\alpha -1} 
      \| \, \partial _x u (\tau) \, \| _{L^{p+1}}^{p+1} % ・ｽ・ｽ ・ｽ・ｽ
      \, \mathrm{d}\tau \\
 & \qquad \quad \quad \; \; \; \: \, 
  + p \int ^t_0 (1+\tau )^\alpha 
    \int _{\partial _x u < 0} 
    f''(u) \left|\, \partial _x u \, \right|^{p+2} 
    \, \mathrm{d}x \mathrm{d}\tau  
       \quad \bigl( t \ge 0 \bigr). 
\end{aligned}
\end{align}
}

\medskip

\noindent
{\bf Proof of Lemma 5.4.}\quad
Multiplying the equation in the problem (1.1), that is, 
$$
\partial_tu +\partial_x \bigl(f(u) \bigr)
  = \mu \, 
    \partial_x \left( \, 
    \left| \, \partial_xu \, \right|^{p-1} \partial_xu \, 
    \right)
$$
by 
$$
- \partial_x \left( \, 
\left| \, \partial_xu \, \right|^{p-1} \partial_xu \, 
\right), 
$$ 
we obtain the divergence form 
\begin{align}
\begin{aligned}
&\partial_t 
 \left(\frac{1}{p+1} \left|\, \partial _x u \, \right|^{p+1} \right) 
 + \partial _x \Bigl( \, 
   - \, \left|\, \partial _x u \, \right|^{p-1} 
   \partial _x u \cdot \partial _t u \, \Bigr) \\
&\qquad \; \, 
 + \partial _x \left( 
 - \, \frac{p}{p+1} \, 
 f'(u) \left|\, \partial _x u \, \right|^{p+1} 
 \, \right) \\
&\qquad \; \, + \frac{p}{p+1} \, 
          f''(u) \left|\, \partial _x u \, \right|^{p+1}\partial _x u 
+\mu \, p\, q \left|\, \partial _x u \, \right|^{2(p-1)} 
   \bigl( \, \partial _x^2 u \, \bigr)^2 
= 0. 
\end{aligned}
\end{align}
Integrating the divergence form (5.14) with respect to $x$, 
we have 
\begin{align}
\begin{aligned}
&\frac{1}{p+1} \, 
\frac{\mathrm{d}}{\mathrm{d}t}\, 
\,\Vert \, \partial _x u (t) \, \Vert_{L^{p+1}}^{p+1} 
+\mu \, p^2 \int ^{\infty }_{-\infty } 
\left|\, \partial _x u \, \right|^{2(p-1)} 
\bigl( \, \partial _x^2 u \, \bigr)^2 
\, \mathrm{d}x \\
&\qquad \qquad \qquad \qquad \quad \; \; \; \: 
 +\frac{p}{p+1} \int ^{\infty }_{-\infty } 
  f''(u) \left|\, \partial _x u \, \right|^{p+1}\partial _x u 
  \, \mathrm{d}x 
= 0. 
\end{aligned}
\end{align}
We separate the integral region 
to the third term on the left-hand side of (5.15) as 
\begin{align}
\begin{aligned}
&\int ^{\infty }_{-\infty } 
 f''(u) \left|\, \partial _x u \, \right|^{p+1}\partial _x u 
 \, \mathrm{d}x \\
&= \int _{\partial _x u \geq 0 } + \int _{\partial _x u < 0 } \\
&= \int _{\partial _x u \geq 0 } 
   f''(u) \left|\, \partial _x u \, \right|^{p+2} \, \mathrm{d}x 
  - \int _{\partial _x u < 0 } 
  f''(u) \left|\, \partial _x u \, \right|^{p+2} \, \mathrm{d}x. 
\end{aligned}
\end{align}
Substituting (5.16) into (5.15), we get the following equality 
\begin{align}
\begin{aligned}
&\frac{1}{p+1} \, 
\frac{\mathrm{d}}{\mathrm{d}t}\, 
\,\Vert \, \partial _x u (t) \, \Vert_{L^{p+1}}^{p+1} 
+\mu \, p^2 \int ^{\infty }_{-\infty } 
\left|\, \partial _x u \, \right|^{2(p-1)} 
\bigl( \, \partial _x^2 u \, \bigr)^2 
\, \mathrm{d}x \\
& \; \: \, +\frac{p}{p+1} \int _{\partial _x u \geq 0 } 
   f''(u) \left|\, \partial _x u \, \right|^{p+2} \, \mathrm{d}x 
= \frac{p}{p+1} \int _{\partial _x u < 0 } 
  f''(u) \left|\, \partial _x u \, \right|^{p+2} \, \mathrm{d}x. 
\end{aligned}
\end{align}
Multiplying (5.17) by 
$(1+t)^{\alpha}$ with $\alpha>0$ 
and integrating over $(0,t)$ with respect to the time, 
we complete the proof of Lemma 5.4. 

\bigskip

\noindent
{\bf Proof of Proposition 5.3.}\quad
We use the following important results (cf. \cite{yoshida'}). 

\medskip

\noindent
{\bf Lemma 5.5.}\quad {\it
For any $s \geq 0$, there exists a positive constant $C_{s}$ 
such that 
\begin{align}
\begin{aligned}
\int _{\partial _x u < 0 } 
f''(u) \left|\, \partial _x u \, \right|^{s} \, \mathrm{d}x 
\leq C_{s}
     \int _{\partial _x u < 0 } 
     \left|\, \partial _x \phi \, \right|^{s} \, \mathrm{d}x. 
\end{aligned}
\end{align}
}

\noindent
In fact, taking care of the relation by using Lemma 2.2 
\begin{align}
\begin{aligned}
\partial _x u = \partial _x U^r + \partial _x \phi <0 \, 
\Longleftrightarrow  \, \partial _x \phi <0, 
                     \, \partial _x U^r < \bigl|\, \partial _x \phi \, \bigr|, 
\end{aligned}
\end{align}
we immediately have 
\begin{align}
\begin{aligned}
\int _{\partial _x u < 0 } 
&f''(u) \left|\, \partial _x u \, \right|^{s} \, \mathrm{d}x \\
&\leq 2^{s} \, 
     \left( \, \max_{ | u | \leq \widetilde{C}} f''(u) \, \right)
     \int _{\partial _x \phi < 0, \partial _x U^r < | \partial _x \phi | } 
     \left|\, \partial _x \phi \, \right|^{s} \, \mathrm{d}x.
\end{aligned}
\end{align}

\medskip
Since $ \partial _x u $ is absolutely continuous, 
we first note that 
for any 
$
x \in \bigl\{ \, x \in \mathbb{R} \, \, \bigr. 
 \bigl| \, \partial _x u \, < \, 0 \, \bigr\}, 
$
there exsists 
$
x_{k} \in \mathbb{R}\cup \{ - \infty \}
$
such that 
$$
\partial _x u (x_{k}) = 0, \; \; 
  \partial _x u (y) \, < \, 0 \; 
  \bigl( \, y \in 
  ( x_{k},x ) \, \bigr).
$$
Therefore, 
it follows that for such $x$ and $x_{k}$ with $q\geq p\, (\, >1\, )$, 
\begin{align}
\left|\, \partial _x u \, \right|^q 
= \left(\, - \partial _x u \, \right)^q 
= q \int_{x_k}^{x} 
   \left(\, - \partial _x u \, \right)^{q-1} \left(\, - \partial _x^2 u \, \right) 
   \, \mathrm{d}y 
\end{align}
By using the Cauchy-Schwarz inequality, we have 

\medskip

\noindent
{\bf Lemma 5.6.}\quad {\it
It holds that 
\begin{align}
\begin{aligned}
& \int _{\partial _x u < 0 } 
  \left|\, \partial _x u \, \right|^{p+2} \, \mathrm{d}x \\
& \le C_{p} 
     \left( \, 
     \int _{\partial _x u < 0 } 
     \left|\, \partial _x u \, \right|^{2(p-1)} 
     \bigl( \, \partial _x^2 u \, \bigr)^2 
     \, \mathrm{d}x \, 
     \right)^{\frac{1}{3p+1}} 
     \left( \, 
     \int _{\partial _x u < 0 } 
     \left|\, \partial _x u \, \right|^{p+1} 
     \, \mathrm{d}x \, 
     \right)^{\frac{3p+2}{3p+1}}. 
\end{aligned}
\end{align}
}
By using Young's inequality to (5.22), we also have 

\medskip

\noindent
{\bf Lemma 5.7.}\quad {\it
It follows that for any $\epsilon>0$, 
there exists a positive constant $C_p({\epsilon})$ such that, 
\begin{align}
\begin{aligned}
& \int _{\partial _x u < 0 } 
  \left|\, \partial _x u \, \right|^{p+2} \, \mathrm{d}x \\
& \le \epsilon 
      \int _{\partial _x u < 0 } 
      \left|\, \partial _x u \, \right|^{2(p-1)} 
      \bigl( \, \partial _x^2 u \, \bigr)^2 
      \, \mathrm{d}x \, 
      + C_p({\epsilon}) \left( \, 
        \int _{\partial _x u < 0 } 
        \left|\, \partial _x u \, \right|^{p+1} 
        \, \mathrm{d}x \, 
        \right)^{\frac{3p+2}{3p}}. 
\end{aligned}
\end{align}
}
By using Lemma 5.5, Lemma 5.6 and Lemma 5.7 
with $\epsilon = \frac{\mu \, p^2 \, (p+1)}{2}$, we have 
\begin{align}
\begin{aligned}
 &(1+t)^\alpha  
  \| \, \partial _x u (t)\, \| _{L^{p+1}}^{p+1} \\
 & \quad 
  + \frac{\mu \, p^2 \, (p+1)}{2} 
    \int ^t_0 (1+\tau )^\alpha  \int _{-\infty}^{\infty} 
    \bigl| \, \partial _x u \, \bigr|^{2(p-1)} 
    \left( \, \partial _x^2 u \, \right)^2 
    \, \mathrm{d}x \mathrm{d}\tau \\
 & \qquad \quad \quad \; \; \; \: \, 
  + p \int ^t_0 (1+\tau )^\alpha  
    \int _{\partial _x u \geq 0} 
    f''(u) \left|\, \partial _x u \, \right|^{p+1} 
    \, \mathrm{d}x \mathrm{d}\tau  \\
 &\leq \| \, \partial _x u_{0} \, \| _{L^{p+1}}^{p+1} \\ % ・ｽ・ｽ ・ｽ・ｽ
 & \quad 
       + \alpha \int ^t_0 (1+\tau )^{\alpha -1} 
         \Bigl( \, 
         \| \, \partial _x \phi (\tau) \, \| _{L^{p+1}}^{p+1} 
         + \| \, \partial _x U^r (\tau) \, \| _{L^{p+1}}^{p+1} 
         \, \Bigr) 
         \, \mathrm{d}\tau \\ 
 & \quad 
  + C_{p} \int ^t_0 (1+\tau )^\alpha 
    \left( \, 
  \int _{\partial _x u < 0 } 
  \left|\, \partial _x u \, \right|^{p+1} 
  \, \mathrm{d}x \, 
  \right)^{\frac{2}{3p}+1} \mathrm{d}\tau . 
\end{aligned}
\end{align}
By using Proposition 5.2, we get the following time-decay estimates. 

\medskip

\noindent
{\bf Lemma 5.8.}\quad {\it
There exist positive constants $\alpha \gg 1$ and $C_{\alpha,p,q}$, 
such that 
 \begin{align}
 \begin{aligned}
 &\int ^t_0 (1+\tau )^\alpha  \int _{-\infty}^{\infty} 
    | \, \phi \, |^{q-2} 
    \bigl( \, \partial _x \phi \, \bigr)^2 
    \left( \, 
    \bigl| \partial _x \phi \bigr|^{p-1} 
    + \bigl| \partial _x U^r \bigr|^{p-1}  
    \, \right) \, \mathrm{d}x \mathrm{d}\tau  \\
 &\leq C\left( \, \alpha, \, p, \, q, \, \phi_0 \, \right) \, 
       (1+t)^{\alpha - \frac{q-1}{2p}}
       \quad \bigl( t \ge 0 \bigr). 
 \end{aligned}
 \end{align}
}

\medskip

\noindent
By using Lemma 5.8 
with $\alpha \mapsto \alpha - 1 \gg 1$ and $q=2$, we have 
\begin{align}
\alpha \int ^t_0 (1+\tau )^{\alpha - 1} 
\| \, \partial _x \phi (\tau) \, \| _{L^{p+1}}^{p+1} 
\, \mathrm{d}\tau  
\leq C\left( \, \alpha, \, p, \, \phi_0 \, \right) \, 
     (1+t)^{\alpha - \frac{2p+1}{2p}}. 
\end{align}
We can also estimate by using Lemma 2.2 as 
\begin{align}
\alpha \int ^t_0 (1+\tau )^{\alpha - 1} 
\| \, \partial _x U^r (\tau) \, \| _{L^{p+1}}^{p+1} 
\, \mathrm{d}\tau  
\leq C\left( \, \alpha, \, p \, \right) \, 
     (1+t)^{\alpha - p}. 
\end{align}
By using the uniform boundedness in Lemma 3.2, that is, 
$$
\| \, \partial _x u (t) 
\, \| _{L^{p+1}}^{p+1} 
%&+\int _0^{\infty} \int _{-\infty}^{\infty} 
% \bigl| \, \partial _x u \, \bigr|^{2(p-1)} 
% \left( \, \partial _x^2 u \, \right)^2 
% \, \mathrm{d}x \mathrm{d}t \\
%&\qquad \quad 
\leq 
 C_{p}\bigl(\, \| \, \phi _0 \, \| _{L^2}, 
 \| \, \partial _x u_0 \, \| _{L^{p+1}} \, \bigr) 
$$
and Lemma 5.8 
with $q=2$, we can estimate as 
\begin{align}
\begin{aligned}
&C_{p} \int ^t_0 (1+\tau )^\alpha 
 \left( \, 
 \int _{\partial _x u < 0 } 
 \left|\, \partial _x u \, \right|^{p+1} 
 \, \mathrm{d}x \, 
 \right)^{\frac{2}{3p}+1} \mathrm{d}\tau \\
&\leq C_{p} \int ^t_0 (1+\tau )^\alpha 
      \int _{\partial _x u < 0 } 
      \left|\, \partial _x \phi \, \right|^{p+1} 
      \, \mathrm{d}x 
      \cdot 
      \| \, \partial _x u (\tau) \, \| _{L^{p+1}}^{\frac{2(p+1)}{3p}} 
      \, \mathrm{d}\tau \\
&\leq C\left( \, p, \, \phi_0, \, \partial _x u_0 \, \right) \, 
      \int ^t_0 (1+\tau )^\alpha 
      \int _{-\infty}^{\infty} 
      \left|\, \partial _x \phi \, \right|^{p+1} 
      \, \mathrm{d}x \mathrm{d}\tau \\
&\leq C\left( \, \alpha, \, p, \, \phi_0, \, \partial _x u_0 \, \right) \,  % ・ｽ・ｽ ・ｽ・ｽ
      (1+t)^{\alpha - \frac{1}{2p}}. 
\end{aligned}
\end{align}
Substituting (5.26), (5.27) and (5.28) into (5.24), 
we complete the proof of Proposition 5.3. 
In particular, we have 
\begin{align}
\| \, \partial _x u (t) \, \| _{L^{p+1}}^{p+1} 
\, \mathrm{d}\tau  
\leq C\left( \, p, \, \phi_0, \, \partial _x u_0 \, \right) \, 
     (1+t)^{- \frac{1}{2p}}, 
\end{align}
and 
\begin{align}
\begin{aligned}
\| \, \partial _x \phi (t) \, \| _{L^{p+1}}^{p+1} 
\, \mathrm{d}\tau  
& \leq \| \, \partial _x u (t) \, \| _{L^{p+1}}^{p+1} 
       + \| \, \partial _x U^r (t) \, \| _{L^{p+1}}^{p+1} \\
& \leq C\left( \, p, \, \phi_0, \, \partial _x u_0 \, \right) \, 
       (1+t)^{- \frac{1}{2p}} 
\end{aligned}
\end{align}
for $1 \leq q < \infty$.

\medskip

{\bf Proof of Theorem 3.2.}\ 
We already have proved the decay estimate 
of $\| \,\phi(t) \,\|_{L^q}$ with $1 \le q < \infty$. 
Therefore we only 
show 
the following time-decay estimate for the higher order derivative
\begin{align}
\begin{aligned}
&\bigl|\bigl|\,
 \partial _x u(t) \,
 \bigr|\bigr|_{L^{p+1} }, \, \; \; 
 \left|\left|\,
 \partial _x \phi (t) 
 \,\right|\right|_{L^{p+1} } \\[5pt]
& \leq \left\{\begin{array} {ll} 
      C( \, \epsilon, \, p, \, \phi_0, \, \partial _x u_0 \, ) \, 
      (1+t)^{-\frac{p}{p+1}}\\[15pt] 
      \,  \, \, \: \; \; \quad \qquad 
      \left( \, 1 < p \le 
      \displaystyle{\frac{1}{3} + 
      \sqrt{\frac{11}{18} - \frac{(p+1)(3p-2)}{3}\, \epsilon } }
      \, \right),\\[15pt] 
      C( \, \epsilon, \, p, \, \phi_0, \, \partial _x u_0 \, ) \, 
      (1+t)^{-\frac{3}{2(p+1)(3p-2)} + \epsilon}\\[15pt] 
      \,  \, \, \; \; \quad \qquad \qquad 
      \left( \,  
      \displaystyle{\frac{1}{3} + 
      \sqrt{\frac{11}{18} - \frac{(p+1)(3p-2)}{3}\, \epsilon } }
      < p \, \right)
      \end{array}
      \right.\,
\end{aligned}
\end{align}
for any $0<\epsilon \ll 1$, 
and the $L^{\infty}$-estimate for $\phi$. 

\medskip

We first prove (5.31). 
Substituting (5.29) into (5.28), 
we have 
\begin{align}
\begin{aligned}
&C_{p} \int ^t_0 (1+\tau )^\alpha 
      \int _{\partial _x u < 0 } 
      \left|\, \partial _x \phi \, \right|^{p+1} 
      \, \mathrm{d}x 
      \cdot 
      \| \, \partial _x u (\tau) \, \| _{L^{p+1}}^{\frac{2(p+1)}{3p}} 
      \, \mathrm{d}\tau \\
&\leq C\left( \, p, \, \phi_0, \, \partial _x u_0 \, \right) \, 
      \int ^t_0 (1+\tau )
      ^{\alpha - \frac{1}{2p} \cdot \frac{2}{3p}}
      \int _{-\infty}^{\infty} 
      \left|\, \partial _x \phi \, \right|^{p+1} 
      \, \mathrm{d}x \mathrm{d}\tau. 
\end{aligned}
\end{align}

\medskip

\noindent
By using Lemma 5.8 
with $\alpha \mapsto \alpha - \frac{1}{2p} \cdot \frac{2}{3p} \gg 1$ 
and $q=2$, we also have 
\begin{align}
\alpha \int ^t_0 (1+\tau )^{\alpha - 1} 
\| \, \partial _x \phi (\tau) \, \| _{L^{p+1}}^{p+1} 
\, \mathrm{d}\tau  
\leq C\left( \, \alpha, \, p, \, \phi_0 \, \right) \, 
     (1+t)^{\alpha - \frac{1}{2p} \cdot \frac{2}{3p} - \frac{1}{2p}}. 
\end{align}
Substituting (5.33) into (5.24), 
we have 
\begin{align}
\begin{aligned}
&\bigl|\bigl|\,
 \partial _x u(t) \,
 \bigr|\bigr|_{L^{p+1} }^{p+1}, \, \; \; 
 \left|\left|\,
 \partial _x \phi (t) 
 \,\right|\right|_{L^{p+1} }^{p+1}  \\[5pt] % ・ｽ・ｽ ・ｽ・ｽ
& \leq C( \, p, \, \phi_0, \, \partial _x u_0 \, ) \\
& \quad \times 
       \left( \, 
       (1+t)^{-\frac{2p+1}{2p}} + (1+t)^{-p} 
       + (1+t)
       ^{-\left( \frac{1}{2p} \cdot \frac{2}{3p} + \frac{1}{2p} \right)} 
       \,  \right) \\[5pt] 
& \leq \left\{\begin{array} {ll} 
      C( \, p, \, \phi_0, \, \partial _x u_0 \, ) \, 
      \left( \, 
       (1+t)^{-p} 
       + (1+t)
       ^{-\left( \frac{1}{2p} \cdot \frac{2}{3p} + \frac{1}{2p} \right)} 
       \,  \right) \\[15pt] 
      \, \: \: \; \; \quad \qquad \qquad \qquad \qquad \qquad \qquad \, 
      \left( \, 1 < p \le \displaystyle{\frac{1 + \sqrt{3}}{2}} \, \right),\\[15pt] 
      C( \, p, \, \phi_0, \, \partial _x u_0 \, ) \, 
      \left( \, 
       (1+t)^{-\frac{2p+1}{2p}} 
       + (1+t)
       ^{-\left( \frac{1}{2p} \cdot \frac{2}{3p} + \frac{1}{2p} \right)} 
       \,  \right) \\[15pt] 
      \, \: \: \; \; \quad \qquad \qquad \qquad \qquad \qquad \qquad \qquad 
      \left( \,  \displaystyle{\frac{1 + \sqrt{3}}{2}} < p \, \right).
      \end{array}
      \right.\,
\end{aligned}
\end{align}
Iterating ``$\infty$''-times the above process, we will get 
\begin{align}
\begin{aligned}
&\bigl|\bigl|\,
 \partial _x u(t) \,
 \bigr|\bigr|_{L^{p+1} }^{p+1}, \, \; \; 
 \left|\left|\,
 \partial _x \phi (t) 
 \,\right|\right|_{L^{p+1} }^{p+1} \\[5pt]
& \leq \left\{\begin{array} {ll} 
      C( \, \epsilon, \, p, \, \phi_0, \, \partial _x u_0 \, ) \, 
      \left( \, 
       (1+t)^{-p} 
       + (1+t)
         ^{- \mathlarger{ \frac{1}{2p}} 
         \substack{{\infty }\\{\substack{{\lsum }\\{n=0}}}}
         \mathlarger{\left( \frac{2}{3p} \right)^{n} + \epsilon} }
       \,  \right) \\[15pt] 
      \quad \quad \quad \; \: \, \, 
      \, \: \: \; \; \quad \qquad \qquad \qquad \qquad \qquad \qquad \, 
      \left( \, 1 < p \le \displaystyle{\frac{1 + \sqrt{3}}{2}} \, \right),\\[15pt] 
      C( \, \epsilon, \, p, \, \phi_0, \, \partial _x u_0 \, ) \, 
      \left( \, 
       (1+t)^{-\frac{2p+1}{2p}} 
       + (1+t)
       ^{- \mathlarger{ \frac{1}{2p}} 
       \substack{{\infty }\\{\substack{{\lsum }\\{n=0}}}}
       \mathlarger{\left(  \frac{2}{3p} \right)^{n}+ \epsilon} }
       \,  \right) \\[15pt] 
      \quad \quad \quad \; \: \, \, 
      \, \: \: \; \; \quad \qquad \qquad \qquad \qquad \qquad \qquad \qquad 
      \left( \,  \displaystyle{\frac{1 + \sqrt{3}}{2}} < p \, \right),
      \end{array}
      \right.\, \\[15pt] 
& \leq \left\{\begin{array} {ll} 
      C( \, \epsilon, \, p, \, \phi_0, \, \partial _x u_0 \, ) \, 
      \left( \, 
       (1+t)^{-p} 
       + (1+t)
       ^{- \frac{1}{2p} \cdot \frac{3p}{3p-2} + \epsilon } 
       \,  \right) \\[15pt] 
       \; \; \: \, \, \, 
      \, \: \: \; \; \quad \qquad \qquad \qquad \qquad \qquad \qquad \, 
      \left( \, 1 < p \le \displaystyle{\frac{1 + \sqrt{3}}{2}} \, \right),\\[15pt] 
      C( \, \epsilon, \, p, \, \phi_0, \, \partial _x u_0 \, ) \, 
      \left( \, 
       (1+t)^{-\frac{2p+1}{2p}} 
       + (1+t)
       ^{- \frac{1}{2p} \cdot \frac{3p}{3p-2} + \epsilon } 
       \,  \right) \\[15pt] 
       \; \; \: \, \, \, 
      \, \: \: \; \; \quad \qquad \qquad \qquad \qquad \qquad \qquad \qquad 
      \left( \,  \displaystyle{\frac{1 + \sqrt{3}}{2}} < p \, \right),
      \end{array}
      \right.\, \\[15pt] 
& \leq \left\{\begin{array} {ll} 
      C( \, \epsilon, \, p, \, \phi_0, \, \partial _x u_0 \, ) \, 
      (1+t)^{-{p}} \, \: \; \qquad 
      \left( \, 1 < p 
      \le \displaystyle{\frac{1}{3} + 
      \sqrt{\frac{11}{18} - \frac{3p-2}{3}\, \epsilon } } \, \right),\\[15pt] 
      C( \, \epsilon, \, p, \, \phi_0, \, \partial _x u_0 \, ) \, 
      (1+t)^{-\frac{3}{2(3p-2)} + \epsilon} \, \: \quad 
      \left( \,  
      \displaystyle{\frac{1}{3} + 
      \sqrt{\frac{11}{18} - \frac{3p-2}{3}\, \epsilon } } 
      < p \, \right)
      \end{array}
      \right.\,
\end{aligned}
\end{align}
for any $0<\epsilon \ll 1$. 
%We also note the following: 
%if $1 < p \le \frac{1 + \sqrt{3}}{2}$, 
%then 
%$$
%(1+t)^{-\frac{2p+1}{2p}}  
%\leq (1+t)^{- p}, 
%$$
%if $p > \frac{1 + \sqrt{3}}{2}$, 
%then 
%$$
%(1+t)^{- p} 
%\leq (1+t)^{-\frac{2p+1}{2p}}
%\leq (1+t)^{-\frac{3}{3p-2}}, 
%$$
%if $p > \frac{2 + \sqrt{22}}{6}$, 
%then 
%$$
%(1+t)^{- p} 
%\leq (1+t)^{-\frac{3}{3p-2}}, 
%$$
%and if $1 < p \leq \frac{2 + \sqrt{22}}{6}$, 
%then 
%$$
%(1+t)^{-\frac{3}{3p-2}} 
%\leq (1+t)^{- p}. 
%$$

\medskip

\noindent
Thus, we get (5.31). 

We finally show the $L^{\infty}$-estimate for $\phi$ 
by using the Gagliardo-Nirenberg inequality.  
Substituting (5.11) and (5.31) into (4.19), we get 
\begin{align}
\begin{aligned}
&\bigl|\bigl|\,
 \phi(t) \,
 \bigr|\bigr|_{L^{\infty } } \\[5pt] 
&\leq \left\{\begin{array} {ll} 
      C( \, \epsilon, \, p, \, \theta, \, \phi_0, \, \partial _x u_0 \, ) \, 
      (1+t)
      ^{-\frac{1}{2p} + 
      \left( \frac{2p+1}{2p(p+1)} - \frac{p}{p+1} \right) \theta} \\[15pt] 
      \, \: \: \; \; \quad \qquad \qquad \qquad \qquad \, 
      \left( \, 1 < p \le 
      \displaystyle{\frac{1}{3} + 
      \sqrt{\frac{11}{18} - \frac{(p+1)(3p-2)}{3}\, \epsilon } }
      \, \right),\\[15pt] 
      C( \, \epsilon, \, p, \, \theta, \, \phi_0, \, \partial _x u_0 \, ) \, 
      (1+t)
      ^{-\frac{1}{2p} + 
      \left( \frac{2p+1}{2p(p+1)} - \frac{3}{2(p+1)(3p-2)}  + \epsilon \right) \theta}\\[15pt] 
      \, \: \: \; \; \quad \qquad \qquad \qquad \qquad \qquad 
      \left( \,  
      \displaystyle{\frac{1}{3} + 
      \sqrt{\frac{11}{18} - \frac{(p+1)(3p-2)}{3}\, \epsilon } }
      < p \, \right)
      \end{array}
      \right.\, %\\[5pt] 
%& \leq C( \, p, \, \theta, \, \phi_0, \, \partial _x u_0 \, ) \, 
%       (1+t)
%       ^{-\left( \frac{1}{2p} \cdot \frac{2}{3p} + \frac{1}{2p} \right)} 
\end{aligned}
\end{align}
for $\theta \in (0,1\, ]$ and any $0<\epsilon \ll 1$. 
Consequently, we do complete the proof of Theorem 3.2.
\bigskip 

\noindent
%%%%%%%%%%%%%%%%%%%%%%%%%%%%%%%%%%%%%%%%%%%%%%%%%%%%%%%%%%%%%%%%%%%%%%%%%%%%%
\section{$L^{r+1}$-estimate for the higher order derivative with $r>p$}
In this section, 
we show the time-decay estimates for the higher order derivative 
in the $L^{r+1}$-norm 
with $r>p$, 
in the case where $\phi_0 \in L^1 \cap L^2$ 
with $\partial _x u_0 \in L^{p+1} \cap L^{r+1}$, 
that is, Theorem 3.3. 

\medskip

\noindent
{\bf Proposition 6.1.}\quad {\it
Suppose the same assumptions in Theorem 3.3. 
For any $r>p$, 
there exist positive constants $\alpha$ and $C_{\alpha,p,r}$, 
such that the unique global solution in time $\phi$ of the Cauchy problem {\rm(3.3)} 
satisfies the following $L^{r+1}$-energy estimate 
\begin{align}
\begin{aligned}
 &(1+t)^\alpha  
  \| \, \partial _x u (t)\, \| _{L^{r+1}}^{r+1} 
  + \int ^t_0 (1+\tau )^\alpha  \int _{-\infty}^{\infty} 
    \bigl| \, \partial _x u \, \bigr|^{p+r-2} 
    \left( \, \partial _x^2 u \, \right)^2 
    \, \mathrm{d}x \mathrm{d}\tau \\
 & \qquad \qquad \qquad \quad \quad \quad \, 
  + \int ^t_0 (1+\tau )^\alpha  
    \| \, \partial _x u (\tau) \, \| _{L^{r+2}}^{r+2} 
    \, \mathrm{d}\tau  \\
 &\leq C_{\alpha,p,r} 
       \| \, \partial _x u_{0} \, \| _{L^{r+1}}^{r+1} \\[5pt] 
 & \; \; 
     + \left\{\begin{array} {ll} 
       C( \, \alpha, \, \epsilon, \, p, \, r, \, \phi_0, \, \partial _x u_0 \, ) \, 
       (1+t)^{\alpha-\frac{2pr+p^2+r}{3p+1}} \\[15pt]  % ・ｽ・ｽ ・ｽ・ｽ・ｽ・ｽ・ｽ・ｽ・ｽ・ｽ・ｽ・ｽ・ｽ・ｽ
       \, \, 
       \left( \, 
       1 < p \le 
       \displaystyle{\frac{1}{3} + 
      \sqrt{\frac{11}{18} - \frac{3p-2}{3}\, \epsilon } }, \; 
       r > p > \displaystyle{\frac{18p^3-17p^2-16p-3}{2(2p+1)}}
       \, \right),\\[15pt] 
       C( \, \alpha, \, \epsilon, \, p, \, r, \, \phi_0, \, \partial _x u_0 \, ) \, 
       (1+t)^{\alpha-\frac{p+2r}{2p(3p-2)} % ・ｽ・ｽ ・ｽ・ｽ・ｽ・ｽ・ｽ・ｽ・ｽ・ｽ・ｽ・ｽ・ｽ・ｽ
       + \frac{2(r-p+1)}{3p} \epsilon} \\[15pt] 
       \, \, \, \: \; \quad \quad \quad \qquad \qquad \qquad \qquad \qquad \qquad 
       \left( \, 
       \displaystyle{\frac{1}{3} + 
       \sqrt{\frac{11}{18} - \frac{3p-2}{3}\, \epsilon } } 
       < p %\le \displaystyle{\frac{2 + \sqrt{7}}{3}}, \; 
%       r > p \ge \displaystyle{\frac{3p^3-2p^2-p}{3p^2-4p-1}}
       \, \right) % ・ｽ・ｽ ・ｽ・ｽ・ｽ・ｽ・ｽ・ｽ・ｽ・ｽ・ｽ・ｽ・ｽ・ｽ
%,\\[15pt] 
%       C( \, \alpha, \, \epsilon, \, p, \, r, \, \phi_0, \, \partial _x u_0 \, ) \, 
%       (1+t)^{\alpha-\frac{6p(r-p)+7p+2r+3}{2(3p+1)(3p-2)(r+1)} 
%       + \frac{p+2r}{2p(3p-2)} \epsilon} \\[15pt] 
%       \, \, \: \; \; \quad \quad \qquad \qquad \qquad \qquad 
%       \left( \, 
%       p > \displaystyle{\frac{2 + \sqrt{7}}{3}}, \; 
%       p < r \le \displaystyle{\frac{3p^3-2p^2-p}{3p^2-4p-1}}
%       \, \right) % ・ｽ・ｽ ・ｽ・ｽ・ｽ・ｽ・ｽ・ｽ・ｽ・ｽ・ｽ・ｽ・ｽ・ｽ
       \end{array}
       \right.\,%\\
%& \quad \qquad \qquad \qquad \qquad \qquad \qquad \qquad \qquad \qquad 
%  \qquad \qquad \qquad 
%  \bigl( t \ge 0 \bigr). 
\end{aligned}
\end{align}
for $t \ge 0$ and any $0<\epsilon \ll 1$. 
}

\medskip

The proof of Proposition 6.1 is given by the following 
three lemmas. 
Because the proofs of them are similar to those of 
Lemma 5.4, Lemma 5.5, Lemma 5.6 and Lemma 5.7, 
we state only here. 

\medskip

\noindent
{\bf Lemma 6.1.}\quad {\it
There exist positive constants 
$C_{p,r}$ and $C_{\alpha,p,r}$ such that 
\begin{align}
\begin{aligned}
 &(1+t)^\alpha  
  \| \, \partial _x u (t)\, \| _{L^{r+1}}^{r+1} \\ % ・ｽ・ｽ ・ｽ・ｽ
 & \quad 
  + \mu \, p\, r \, (r+1) 
    \int ^t_0 (1+\tau )^\alpha  \int _{-\infty}^{\infty} 
    \bigl| \, \partial _x u \, \bigr|^{p+r-2} 
    \left( \, \partial _x^2 u \, \right)^2 
    \, \mathrm{d}x \mathrm{d}\tau \\
 & \qquad \quad \quad \; \; \; \: \, 
  + r \int ^t_0 (1+\tau )^\alpha  
    \int _{\partial _x u \geq 0} 
    f''(u) \left|\, \partial _x u \, \right|^{r+2} 
    \, \mathrm{d}x \mathrm{d}\tau  \\
 &\leq \| \, \partial _x u_{0} \, \| _{L^{r+1}}^{r+1} \\ % ・ｽ・ｽ ・ｽ・ｽ
 & \quad \; \; \; \: \, 
       + C_{\alpha,p,r} \int ^t_0 
         (1+\tau )
         ^{\alpha -{\frac{2p+r+1}{3p+1}}} 
         \left( \, 
         \int _{-\infty}^{\infty} 
         \left|\, \partial _x u \, \right|^{p+1} 
         \, \mathrm{d}x \, 
         \right)^{\frac{p+2r+1}{3p+1}}
         \, \mathrm{d}\tau \\
 & \quad \; \; \; \: \, 
  + C_{p,r} \int ^t_0 (1+\tau )^\alpha 
    \left( \, 
         \int _{\partial _x u < 0 } 
         \left|\, \partial _x u \, \right|^{p+1} 
         \, \mathrm{d}x \, 
         \right)^{\frac{p+2r+2}{3p}}
         \, \mathrm{d}\tau  
       \quad \bigl( t \ge 0 \bigr). 
\end{aligned}
\end{align}
}

\bigskip

\noindent
{\bf Lemma 6.2.}\quad {\it
Assume $p>1$ and $r>p$. 
We have the following interpolation inequalities.  

 \noindent
 {\rm (1)}\ \ There exists a positive 
 constant $C_{p,r}$  such that 
 \begin{align*}
 \begin{aligned}
 &\Vert \, \partial _x u (t) \, \Vert _{L^{r+1} }
 \leq 
  C_{p,r} \left( \, 
     \int _{-\infty}^{\infty}
     \left|\, \partial _x u \, \right|^{p+r-2} % ・ｽ・ｽ ・ｽ・ｽ・ｽ・ｽ・ｽ・ｽ・ｽ・ｽ・ｽ・ｽ・ｽ・ｽ・ｽ・ｽ・ｽ・ｽ・ｽ・ｽ
     \bigl( \, \partial _x^2 u \, \bigr)^2 
     \, \mathrm{d}x \, 
     \right)^{\frac{r-p}{(2p+r+1)(r+1)}} \\
 & \: \; \; \quad \quad \quad \qquad \qquad \times 
     \left( \, 
     \int _{-\infty}^{\infty} % ・ｽ・ｽ ・ｽ・ｽ
     \left|\, \partial _x u \, \right|^{p+1} 
     \, \mathrm{d}x \, 
     \right)^{\frac{p+2r+1}{(2p+r+1)(r+1)}}.
% \quad \bigl( t \ge 0 \bigr). 
 \end{aligned}
 \end{align*}

 \noindent
 {\rm (2)}\ \ 
 There exists a positive 
 constant $C_{p,r}$ 
 such that 
 \begin{align*}
 \begin{aligned}
 &\Vert \, \partial _x u (t) \, \Vert _{L^\infty }
 \leq 
  C_{p,r} \left( \, 
     \int _{-\infty}^{\infty} 
     \left|\, \partial _x u \, \right|^{p+r-2} % ・ｽ・ｽ ・ｽ・ｽ・ｽ・ｽ・ｽ・ｽ・ｽ・ｽ・ｽ・ｽ・ｽ・ｽ・ｽ・ｽ・ｽ・ｽ・ｽ・ｽ
     \bigl( \, \partial _x^2 u \, \bigr)^2 
     \, \mathrm{d}x \, 
     \right)^{\frac{1}{2p+r+1}} \\
 & \: \quad \quad \quad \qquad \qquad \times 
     \left( \, 
     \int _{-\infty}^{\infty} % ・ｽ・ｽ ・ｽ・ｽ
     \left|\, \partial _x u \, \right|^{p+1} 
     \, \mathrm{d}x \, 
     \right)^{\frac{1}{2p+r+1}}.
% \quad \bigl( t \ge 0 \bigr). 
 \end{aligned}
 \end{align*}
}

\bigskip

\noindent
{\bf Lemma 6.3.}\quad {\it
Assume $p>1$ and $r>p$. 
We have the following interpolation inequalities.  

 \noindent
 {\rm (1)}\ \ There exists a positive 
 constant $C_{p,r}$  such that 
 \begin{align*}
 \begin{aligned}
 &\Vert \, \partial _x u (t) \, \Vert 
  _{L^{r+1}_{x} \left( \{ \partial _x u < 0 \} \right) }
 \leq 
  C_{p,r} \left( \, 
     \int _{\partial _x u < 0 } 
     \left|\, \partial _x u \, \right|^{p+r-2} % ・ｽ・ｽ ・ｽ・ｽ・ｽ・ｽ・ｽ・ｽ・ｽ・ｽ・ｽ・ｽ・ｽ・ｽ・ｽ・ｽ・ｽ・ｽ・ｽ・ｽ
     \bigl( \, \partial _x^2 u \, \bigr)^2 
     \, \mathrm{d}x \, 
     \right)^{\frac{r-p}{(2p+r+1)(r+1)}} \\
 & \, \: \: \quad \qquad \qquad \qquad \qquad \qquad \times 
     \left( \, 
     \int _{\partial _x u < 0 } 
     \left|\, \partial _x u \, \right|^{p+1} 
     \, \mathrm{d}x \, 
     \right)^{\frac{p+2r+1}{(2p+r+1)(r+1)}}.
% \quad \bigl( t \ge 0 \bigr). 
 \end{aligned}
 \end{align*}

 \noindent
 {\rm (2)}\ \ 
 There exists a positive 
 constant $C_{p,r}$ 
 such that 
 \begin{align*}
 \begin{aligned}
 &\Vert \, \partial _x u (t) \, \Vert 
  _{L^{\infty}_{x} \left( \{ \partial _x u < 0 \} \right) }
 \leq 
  C_{p,r} \left( \, 
     \int _{\partial _x u < 0 } 
     \left|\, \partial _x u \, \right|^{p+r-2} % ・ｽ・ｽ ・ｽ・ｽ・ｽ・ｽ・ｽ・ｽ・ｽ・ｽ・ｽ・ｽ・ｽ・ｽ・ｽ・ｽ・ｽ・ｽ・ｽ・ｽ
     \bigl( \, \partial _x^2 u \, \bigr)^2 
     \, \mathrm{d}x \, 
     \right)^{\frac{1}{2p+r+1}} \\
 & \, \, \: \: \; \qquad \qquad \qquad \qquad \qquad \times 
     \left( \, 
     \int _{\partial _x u < 0 } 
     \left|\, \partial _x u \, \right|^{p+1} 
     \, \mathrm{d}x \, 
     \right)^{\frac{1}{2p+r+1}}.
% \quad \bigl( t \ge 0 \bigr). 
 \end{aligned}
 \end{align*}
}

\medskip

{\bf Proof of Proposition 6.1.}\ 
By using (5.31), 
we estimate the each terms on the right-hand side of (6.2) as 
\begin{align}
\begin{aligned}
&C_{\alpha,p,r} \int ^t_0 
         (1+\tau )
         ^{\alpha -{\frac{2p+r+1}{3p+1}}} 
         \left( \, 
         \int _{-\infty}^{\infty} 
         \left|\, \partial _x u \, \right|^{p+1} 
         \, \mathrm{d}x \, 
         \right)^{\frac{p+2r+1}{3p+1}}
         \, \mathrm{d}\tau \\[5pt]
& \leq \left\{\begin{array} {ll} 
      C( \, \alpha, \, \epsilon, \, p, \, r, \, \phi_0, \, \partial _x u_0 \, ) \, 
      \displaystyle{ \int ^t_0 
      (1+\tau )
      ^{\alpha - \frac{2p+r+1}{3p+1} -{\frac{p(p+2r+1)}{3p+1}}} 
      \, \mathrm{d}\tau } \\[15pt] 
      \, \: \: \; \; \quad \qquad \qquad \qquad \qquad \qquad \qquad \, 
      \left( \, 1 < p \le 
      \displaystyle{\frac{1}{3} + 
      \sqrt{\frac{11}{18} - \frac{3p-2}{3}\, \epsilon } }
      \, \right),\\[15pt] 
      C( \, \alpha, \, \epsilon, \, p, \, r, \, \phi_0, \, \partial _x u_0 \, ) \, 
      \displaystyle{ \int ^t_0 
      (1+\tau )
      ^{\alpha - \frac{2p+r+1}{3p+1} -{\frac{3(p+2r+1)}{2(3p+1)(3p-2)}}
      + \frac{p+2r+1}{3p+1} \epsilon} 
      \, \mathrm{d}\tau } \\[15pt] 
      \, \: \: \; \; \quad \qquad \qquad \qquad \qquad \qquad \qquad \qquad 
      \left( \,  \displaystyle{\frac{1}{3} + 
      \sqrt{\frac{11}{18} - \frac{3p-2}{3}\, \epsilon } }
      < p \, \right)
      \end{array}
      \right.\, \\[5pt] 
& \leq \left\{\begin{array} {ll} 
      C( \, \alpha, \, \epsilon, \, p, \, r, \, \phi_0, \, \partial _x u_0 \, ) \, 
      (1+t )
      ^{\alpha - \frac{2pr+p^2+r}{3p+1} } \\[15pt] 
      \, \: \: \; \; \quad \qquad \qquad \qquad \qquad \qquad \qquad \, 
      \left( \, 1 < p \le 
      \displaystyle{\frac{1}{3} + 
      \sqrt{\frac{11}{18} - \frac{3p-2}{3}\, \epsilon } }
      \, \right),\\[15pt] 
      C( \, \alpha, \, \epsilon, \, p, \, r, \, \phi_0, \, \partial _x u_0 \, ) \, 
      (1+t )
      ^{\alpha - \frac{6p(r-p)+7p+2r+3}{2(3p+1)(3p-2)} 
      + \frac{p+2r+1}{3p+1} \epsilon} \\[15pt] 
      \, \: \: \; \; \quad \qquad \qquad \qquad \qquad \qquad \qquad \qquad 
      \left( \,  \displaystyle{\frac{1}{3} + 
      \sqrt{\frac{11}{18} - \frac{3p-2}{3}\, \epsilon } }
      < p \, \right),
      \end{array}
      \right.\, %\\[5pt] 
%& \leq \left\{\begin{array} {ll} 
%      C( \, p, \, \phi_0, \, \partial _x u_0 \, ) \, 
%      (1+t)^{-{p}} \, \: \; \quad \qquad 
%      \left( \, 1 < p \le \displaystyle{\frac{2 + \sqrt{22}}{6}} \, \right),\\[10pt] 
%      C( \, p, \, \phi_0, \, \partial _x u_0 \, ) \, 
%      (1+t)^{-\frac{3}{2(3p-2)}} \, \, \: \quad \quad 
%      \left( \,  \displaystyle{\frac{2 + \sqrt{22}}{6}} < p \, \right).
%      \end{array}
%      \right.\,
\end{aligned}
\end{align}
\begin{align}
\begin{aligned}
&C_{p,r} \int ^t_0 (1+\tau )^\alpha 
    \left( \, 
    \int _{\partial _x u < 0 } 
    \left|\, \partial _x u \, \right|^{p+1} 
    \, \mathrm{d}x \, 
    \right)^{\frac{p+2r+2}{3p}}
    \, \mathrm{d}\tau  \\
& \leq C_{p,r} \int ^t_0 (1+\tau )^\alpha 
       \left( \, 
       \int _{-\infty}^{\infty}
       \left|\, \partial _x \phi \, \right|^{p+1} 
       \, \mathrm{d}x \, 
       \right)
       \| \, \partial _x u (\tau) \, \| _{L^{p+1}}
       ^{\frac{2(p+1)(r-p+1)}{3p}} 
       \, \mathrm{d}\tau \\[5pt]
& \leq \left\{\begin{array} {ll} 
      C( \, \alpha, \, \epsilon, \, p, \, r, \, \phi_0, \, \partial _x u_0 \, ) \, 
      \displaystyle{ \int ^t_0 
      (1+\tau )
      ^{\alpha - \frac{2(r-p+1)}{3} } 
      \| \, \partial _x \phi (\tau) \, \| _{L^{p+1}}^{p+1}
      \, \mathrm{d}\tau } \\[15pt] 
      \, \: \: \; \; \; \; \; \quad \qquad \qquad \qquad \qquad \qquad \qquad \, 
      \left( \, 1 < p \le 
      \displaystyle{\frac{1}{3} + 
      \sqrt{\frac{11}{18} - \frac{3p-2}{3}\, \epsilon } }
      \, \right),\\[15pt] 
      C( \, \alpha, \, \epsilon, \, p, \, r, \, \phi_0, \, \partial _x u_0 \, ) \, 
      \displaystyle{ \int ^t_0 
      (1+\tau )
      ^{\alpha - \frac{r-p+1}{p(3p-2)} + \frac{2(r-p+1)}{3p} \epsilon} 
      \| \, \partial _x \phi (\tau) \, \| _{L^{p+1}}^{p+1}
      \, \mathrm{d}\tau } \\[15pt] 
      \, \: \: \; \; \; \; \; \quad \qquad \qquad \qquad \qquad \qquad \qquad \qquad 
      \left( \,  \displaystyle{\frac{1}{3} + 
      \sqrt{\frac{11}{18} - \frac{3p-2}{3}\, \epsilon } }
      < p \, \right)
      \end{array}
      \right.\, 
\end{aligned}
\end{align}
for any $0<\epsilon \ll 1$. \\
\noindent
By using Lemma 5.8 with 
\begin{align*}
\begin{aligned}
\alpha \mapsto \left\{\begin{array} {ll} 
               \alpha-\displaystyle{\frac{2(r-p+1)}{3}} 
               \, \, \: \; \quad \quad \qquad \qquad 
               \left( \, 1 < p \le 
               \displaystyle{\frac{1}{3} + 
               \sqrt{\frac{11}{18} - \frac{3p-2}{3}\, \epsilon } } 
               \, \right),\\[15pt] 
               \alpha
               -\displaystyle{\left( \, 
                \frac{r-p+1}{p(3p-2)} - \frac{2(r-p+1)}{3p} \epsilon 
                \, \right)}
               \quad 
               \left( \,  \displaystyle{\frac{1}{3} + 
               \sqrt{\frac{11}{18} - \frac{3p-2}{3}\, \epsilon } }
               < p \, \right)
               \end{array}
               \right.\, 
\end{aligned}
\end{align*}
and $q=2$, 
we get 
\begin{align}
\begin{aligned}
&C_{p,r} \int ^t_0 (1+\tau )^\alpha 
    \left( \, 
    \int _{\partial _x u < 0 } 
    \left|\, \partial _x u \, \right|^{p+1} 
    \, \mathrm{d}x \, 
    \right)^{\frac{p+2r+2}{3p}}
    \, \mathrm{d}\tau  \\[5pt]
& \leq \left\{\begin{array} {ll} 
      C( \, \alpha, \, \epsilon, \, p, \, r, \, \phi_0, \, \partial _x u_0 \, ) \, 
      (1+t )
      ^{\alpha - \frac{4p(r-p)+4p+3}{6p} } \\[15pt] 
      \, \: \: \; \; \quad \qquad \qquad \qquad \qquad \, 
      \left( \, 1 < p \le 
      \displaystyle{\frac{1}{3} + 
      \sqrt{\frac{11}{18} - \frac{3p-2}{3}\, \epsilon } } 
      \, \right),\\[15pt] 
      C( \, \alpha, \, \epsilon, \, p, \, r, \, \phi_0, \, \partial _x u_0 \, ) \, 
      (1+t )
      ^{\alpha - \frac{p+2r}{2p(3p-2)} + \frac{2(r-p+1)}{3p} \epsilon} \\[15pt] 
      \, \: \: \; \; \quad \qquad \qquad \qquad \qquad \qquad 
      \left( \,  \displaystyle{\frac{1}{3} + 
      \sqrt{\frac{11}{18} - \frac{3p-2}{3}\, \epsilon } }
      < p \, \right)
      \end{array}
      \right.\, 
\end{aligned}
\end{align}
for any $0<\epsilon \ll 1$. \\
\noindent
Substituting (6.3) and (6.5) into (6.2), 
we have 
\begin{align}
\begin{aligned}
 &(1+t)^\alpha  
  \| \, \partial _x u (t)\, \| _{L^{r+1}}^{r+1} \\ % ・ｽ・ｽ ・ｽ・ｽ・ｽ・ｽ
 & \quad 
  + \int ^t_0 (1+\tau )^\alpha  \int _{-\infty}^{\infty} 
    \bigl| \, \partial _x u \, \bigr|^{p+r-2} 
    \left( \, \partial _x^2 u \, \right)^2 
    \, \mathrm{d}x \mathrm{d}\tau \\
 & \quad 
  + \int ^t_0 (1+\tau )^\alpha  
    \int _{\partial _x u \geq 0} 
    f''(u) \left|\, \partial _x u \, \right|^{r+2} 
    \, \mathrm{d}x \mathrm{d}\tau  \\
& \leq C_{\alpha,p,r} 
       \| \, \partial _x u_{0} \, \| _{L^{r+1}}^{r+1} \\[5pt] % ・ｽ・ｽ ・ｽ・ｽ・ｽ・ｽ
& \quad + 
      \left\{\begin{array} {ll} 
      C( \, \alpha, \, \epsilon, \, p, \, r, \, \phi_0, \, \partial _x u_0 \, ) 
      (1+t)^{\alpha}\\[5pt]
      \times 
      \left( \, 
       (1+t)
       ^{- \frac{2pr+p^2+r}{3p+1}} 
       + (1+t)
       ^{- \frac{4p(r-p)+4p+3}{6p} } 
       \,  \right) \\[15pt] 
      \, \: \: \; \; \quad \qquad \qquad \qquad \qquad \qquad \qquad \, 
      \left( \, 1 < p \le 
      \displaystyle{\frac{1}{3} + 
      \sqrt{\frac{11}{18} - \frac{3p-2}{3}\, \epsilon } }
      \, \right),\\[15pt] 
      C( \, \alpha, \, \epsilon, \, p, \, r, \, \phi_0, \, \partial _x u_0 \, ) 
      (1+t)^{\alpha}\\[5pt]
      \times 
      \left( \, 
       (1+t)
       ^{- \frac{6p(r-p)+7p+2r+3}{2(3p+1)(3p-2)} 
       + \frac{p+2r+1}{3p+1} \epsilon} 
       + (1+t)
       ^{- \frac{p+2r}{2p(3p-2)} 
       + \frac{2(r-p+1)}{3p} \epsilon} 
       \,  \right) \\[15pt] 
      \, \: \: \; \; \quad \qquad \qquad \qquad \qquad \qquad \qquad \qquad 
      \left( \,  \displaystyle{\frac{1}{3} + 
      \sqrt{\frac{11}{18} - \frac{3p-2}{3}\, \epsilon } }
      < p \, \right)
      \end{array}
      \right.\, \\[5pt] 
\end{aligned}
\end{align}
for any $0<\epsilon \ll 1$.\\
\noindent
We also note the following: 
if $1 < p 
    \le \frac{1}{3} + 
        \sqrt{\frac{11}{18} - \frac{3p-2}{3}\, \epsilon} $, 
then $r > p > \frac{18p^3-17p^2-16p-3}{2(2p+1)}$ and 
$$
(1+t)^{- \frac{4p(r-p)+4p+3}{6p} } 
\leq (1+t)^{- \frac{2pr+p^2+r}{3p+1}},  % ・ｽ・ｽ ・ｽ・ｽ・ｽ・ｽ・ｽ・ｽ・ｽ・ｽ・ｽ・ｽ・ｽ・ｽ
$$
and if $\epsilon  
      < \frac{3(r-p)(p-1)(3p+1)}{(3p-2)\, \left| \, 9p^2-p-2r-2 \, \right|}$, 
then 
$$
(1+t)^{- \frac{6p(r-p)+7p+2r+3}{2(3p+1)(3p-2)} 
+ \frac{p+2r+1}{3p+1} \epsilon} 
\leq (1+t)^{- \frac{p+2r}{2p(3p-2)} 
+ \frac{2(r-p+1)}{3p} \epsilon} \; \; 
( \, \forall p>1, \, \forall r>p  \, ).  % ・ｽ・ｽ ・ｽ・ｽ・ｽ・ｽ・ｽ・ｽ・ｽ・ｽ・ｽ・ｽ・ｽ・ｽ
$$
%and if $p>\frac{2 + \sqrt{7}}{3}$, 
%then $p < r \le \frac{3p^3-2p^2-p}{3p^2-4p-1}$ and 
%$$
%(1+t)^{- \frac{p+2r}{2p(3p-2)} } 
%\leq (1+t)^{- \frac{6p(r-p)+7p+2r+3}{2(3p+1)(3p-2)}
%+ \frac{p+2r+1}{3p+1} \epsilon}. 
%$$

\medskip

\noindent
Thus, we do complete the proof of Proposition 6.1. 

\bigskip 

\noindent
%%%%%%%%%%%%%%%%%%%%%%%%%%%%%%%%%%%%%%%%%%%%%%%%%%%%%%%%%%%
\section{Discussion}
%In this work we studied the time-decay estimates. 
In this section, we discuss the time-decay rates in our main theorems. 
To do that, we recall the time-decay rates of solutions to 
a Cauchy problem for the simplest $p$-Laplacian evolution equation without convestive term: 
\begin{eqnarray}
 \left\{\begin{array}{ll}
  \partial_tu - 
    \mu \, 
    \partial_x \left( \, 
    \left| \, \partial_xu \, \right|^{p-1} \partial_xu \, 
    \right) =0 
  \qquad &(t>0, x\in \mathbb{R}), \\[5pt]
  u(0,x) = u_0(x) \qquad &( x \in \mathbb{R} ),\\[5pt]
  \displaystyle{\lim_{x\to \pm \infty}} u(t,x) =0
  \qquad &\bigl( t \ge 0 \bigr),   
 \end{array}
 \right.\,
\end{eqnarray}
%for $p>1$. 
where, $u=u(t,x)$ denotes the unknown function 
of $t>0$ and $x\in \mathbb{R}$. 
The theorems concerning the time-decay estimates to the problem (7.1) 
are as follows (the proofs are similar to those in the previous sections). 

\medskip

\noindent
{\bf Theorem 7.1.}\quad{\it
If the initial data satisfies 
$u_0 \in L^2$ and 
$\partial _x u_0 \in L^{p+1}$. 
then there uniguely exists a global solution in time $u$ 
of the Cauchy problem {\rm(7.1)} 
satisfying 
\begin{eqnarray*}
\left\{\begin{array}{ll}
u \in C^0\bigl( \, [\, 0,\infty)\, ;L^2 \bigr)
         \cap L^{\infty}\bigl( \, \mathbb{R}^{+} \, ;L^{2} \bigr),\\[5pt]
\partial_x \left( \, 
    \left| \, \partial_xu \, \right|^{p-1} \partial_xu \, 
    \right)
\in L^{2}\bigl(\, {\mathbb{R}^{+}_{t}} \times {\mathbb{R}}_{x} \bigr)
\end{array} 
\right.\,
\end{eqnarray*}
and 
%$$
%\displaystyle{\lim_{t\to \infty}}\sup_{x\in \mathbb{R}}\, 
%\bigl| \, u(t,x)\,\bigr| = 0. 
%$$
%The the solution further 
%satisfies 
the following time-decay estimates 
\begin{eqnarray*}
\left\{\begin{array} {ll}
\left|\left|\,
u (t) 
\, \right|\right|_{L^q }
\leq C( \, p, \, q, \, u_0 \, ) \, 
     (1+t)^{-\frac{1}{3p+1}\left(1-\frac{2}{q}\right)},\\[7pt]
\left|\left|\,
u (t) 
\, \right|\right|_{L^{\infty} }
\leq C( \, \epsilon, \, p, \, q, \, u_0, \, \partial _x u_0 \, ) \, 
     (1+t)^{-\frac{1}{3p+1}+\epsilon} 
\end{array}
  \right.\,
\end{eqnarray*}
for 
$q \in [\, 2, \infty) $ and any $\epsilon>0$. 
}

\medskip

\noindent
{\bf Theorem 7.2.}\quad{\it
If the initial data further satisfies 
$u_0 \in L^1 \cap L^2$ and 
$\partial _x u_0 \in L^{p+1}$, 
then it holds that the unique global solution in time $u$ 
of the Cauchy problem {\rm(7.1)} 
satisfies the following time-decay estimates 
\begin{eqnarray*}
\left\{\begin{array} {ll}
\left|\left| \,
u (t) 
\, \right|\right|_{L^q }
\leq C( \, p, \, q, \, u_0 \, ) \, 
     (1+t)^{-\frac{1}{2p}\left(1-\frac{1}{q}\right)},\\[7pt]
\left|\left|\,
u (t) 
\,\right|\right|_{L^{\infty} }
\leq C( \, p, \, q, \, u_0, \, \partial _x u_0 \, ) \, 
     (1+t)^{-\frac{1}{2p}} % ・ｽ・ｽ ・ｽ・ｽ・ｽ・ｽ・ｽ・ｽ・ｽ・ｽ・ｽ・ｽ・ｽ・ｽNO・ｽ・ｽ
\end{array}
  \right.\,
\end{eqnarray*}
for 
$q \in [\, 1, \infty) $. 
Furthermore, the solution satisfies 
the following time-decay estimates for the higher order derivative
\begin{equation*}
\bigl|\bigl|\,
 \partial _x u(t) \,
 \bigr|\bigr|_{L^{p+1} } 
 \leq C( \, p, \, u_0, \, \partial _x u_0 \, ) \, 
      (1+t)^{-\frac{2p+1}{2p(p+1)}}. 
\end{equation*}
}

\medskip

\noindent
{\bf Theorem 7.3.}\quad{\it
If the initial data further satisfies 
$u_0 \in L^1 \cap L^2$ and 
$\partial _x u_0 \in L^{p+1} \cap L^{r+1} \, (r>p)$, 
then it holds that the unique global solution in time $u$ 
of the Cauchy problem {\rm(7.1)} 
satisfies the following time-decay estimates 
for the higher order derivative
\begin{equation*}
\bigl|\bigl|\,
 \partial _x u(t) \,
 \bigr|\bigr|_{L^{r+1} } 
 \leq C( \, p, \, r, \, u_0, \, \partial _x u_0 \, ) \, 
      (1+t)^{-\frac{6pr+3p+2r+1}{2p(3p+1)(r+1)}}. 
\end{equation*}
}
\medskip

It is clear that the time-decay rates 
in the $L^q$-norm ($2 \leq q \leq \infty$ or $1 \leq q \leq \infty$)
for the lower order $u-\tilde{u}$ or $u-u^r$ in Theorems 1.1, 1.2, 1.4 and 1.5 
are quite or almost the same as $u$ in Theorem 7.1 and 7.2. 
This shows that the affection to the time-decay from the formulation of the equation 
$$
\partial_tu - 
    \mu \, 
    \partial_x \left( \, 
    \left| \, \partial_xu \, \right|^{p-1} \partial_xu \, 
    \right) =0
$$
is stronger than those from 
the asymptotic states, the rarefaction wave $u^r$ (or $U^r$) 
or the constant states $\tilde{u}$ 
(and also from the shape of the flux function $f$, see (4.2), (5.3) and (5.10)). 
In fact, the time-decay in (4.16) is faster than that in (4.15) 
without $\alpha \gg 1$ in Section 4. 
%Therefore it is natural for us to predict those rates are optimal. 
On the other hand, the time-decay rates 
in the $L^{p+1}$-norm or the $L^{r+1}$-norm ($r>p$) 
for the higher order $\partial _x u$ or $\partial _x u-\partial _x u^r$ 
in Theorems 1.2, 1.3, 1.5, 1.6, 7.2 and 7.3 are all different from with each other. 
The reason for the difference must arise from 
that the asumptions for the flux function $f$ and the far field states $u_{\pm}, \, \tilde{u}$, 
and the characteristic propaties of the asymptotic states ($u^r$ or $\tilde{u}$) affect (in some sense)
the strong nonlinearlity 
of the higher order $\partial _x u$ (not $u$), that is, $p$-Laplacian type viscosity 
(see (5.18), (5.19), (5.25), (5.26), (5.33) in Section 5 and (6.5) in Section 6). 
However, the optimality of the all time-decay rates still remains open.

\bigskip

{\bf Acknowledgement.}\quad
%The author would like to thank 
%an anonymous referee for his/her 
%many significant comments and kind advices 
%on our manuscript. 
The auther %also 
thanks 
Professor Akitaka Matsumura for 
his significant comments and kind advices.

% \linenumbers


%%%%%%%%%%%%%%
%\section{}
%\label{}



%% References
%%   \citet{key}         ==>   Author [#]


\bibliographystyle{model6-num-names}
\bibliography{<your-bib-database>}

%% not want to use model6-num-names.bst.





% \end{thebibliography}

\end{document}